\documentclass[reqno,11pt,english]{amsart}
\usepackage[T1]{fontenc}
\usepackage[latin9]{inputenc}
\usepackage{amstext}
\usepackage{amsthm}
\usepackage{amssymb}
\usepackage{comment}

\makeatletter
\usepackage{amsmath,amsfonts,amssymb,amsthm,epsfig}

\usepackage{color}
\usepackage[final,allcolors=blue,colorlinks=true]{hyperref}

\def\R{\mathbb R}
\def\N{\mathbb N}

\def\al{\alpha}
\def\be{\beta}

\def\de{\delta}
\def\ep{\epsilon}
\def\la{\lambda}

\def\om{\omega}
\def\na{\nabla}
\def\Om{\Omega}  
\def\De{\Delta}      

\def\cal{\mathcal}
\def\wq{\infty}
\def\pa{\partial}

\def\loc{\text{\rm loc}}

\def\diam{\text{\rm diam}}

\newcommand{\D}{{\rm d}}

\newcommand{\fint}{-\kern-,400cm\int}

\numberwithin{equation}{section}
\newtheorem{theorem}{Theorem}[section]
\newtheorem*{theorem*}{Theorem}  

\newtheorem{claim}[theorem]{Claim}

\newtheorem*{conclusion*}{Conclusin}

\newtheorem{corollary}[theorem]{Corollary}
\newtheorem*{corollary*}{Corollary}

\newtheorem{lemma}[theorem]{Lemma}
\newtheorem*{lemma*}{Lemma}

\newtheorem*{notation*}{Notation}
\newtheorem{problem}[theorem]{Problem}
\newtheorem{proposition}[theorem]{Proposition}
\newtheorem*{proposition*}{Proposition}

\newtheorem*{remark*}{Remark}

\newtheorem*{example*}{Example}                
\theoremstyle{definition}

\makeatother

\usepackage{babel}
\begin{document}
	\title[$L^{p}$-theory]{$L^p$-regularity  for fourth order elliptic systems
	with antisymmetric potentials in higher dimensions}
	
	\author[C.-Y. Guo, C.Y. Wang and C.-L. Xiang]{Chang-Yu Guo, Changyou Wang and Chang-Lin Xiang$^\ast$}

	
\address[Chang-Yu Guo]{Research Center for Mathematics and Interdisciplinary Sciences, Shandong University 266237,  Qingdao, P. R. China}
\email{changyu.guo@sdu.edu.cn}

\address[Changyou Wang]{Department of Mathematics, Purdue University, West Lafayette, IN 47907, USA}
\email{wang2482@purdue.edu}

\address[Chang-Lin Xiang]{Three Gorges Mathematical Research Center, China Three Gorges University,  443002, Yichang,  P. R. China}
\email{changlin.xiang@ctgu.edu.cn}

\thanks{*Corresponding author: Chang-Lin Xiang}
\thanks{C.-Y. Guo is supported by the Qilu funding of Shandong University (No.~62550089963197), the Natural Science Foundation of Shandong Province (No.~ZR2021QA003) and the National Natural Science Foundation of China (No.~1210010442).
C. Y. Wang is partially supported by NSF grants 1764417 and 2101224.
C.-L. Xiang is financially supported by the National Natural Science Foundation of China (No.~11701045).}
	
	\begin{abstract} We establish an optimal $L^p$-regularity theory for solutions to fourth order
	elliptic systems with antisymmetric potentials in all supercritical dimensions $n\ge 5$:
	$$
	\Delta^2 u=\Delta(D\cdot\nabla u)+div(E\cdot\nabla u)+(\Delta\Omega+G)\cdot\nabla u +f
\qquad	\ {\rm{in}}\ B^n,
	$$
	where $\Omega\in W^{1,2}(B^n, so_m)$ is antisymmetric and $f\in L^p(B^n)$,
	and $D, E, \Omega, G$ satisfy the growth condition (GC-4), under the smallness condition
	of a critical scale invariant norm of $\nabla u$ and $\nabla^2 u$. This system was brought into lights
	from the study of regularity of (stationary) biharmonic maps between manifolds by
	Lamm-Rivi\`ere, Struwe, and Wang.  In particular, our results improve Struwe's H\"older regularity theorem
	to any H\"older exponent $\alpha\in (0,1)$ when $f\equiv 0$, and have applications to both
	approximate biharmonic maps and heat flow of biharmonic maps.

        As a by-product of the techniques, we also extend the $L^p$-regularity theory of harmonic maps by Moser
        to Rivi\`ere-Struwe's second order elliptic systems with antisymmetric
        potentials under the growth condition (GC-2) in all dimensions, which confirms an expectation by Sharp.
	\end{abstract}
	
	\maketitle
	
	{\small
		\smallskip
		\newline
		\tableofcontents}
	\bigskip
	
\section{Introduction and main results}

\subsection{Background and motivation}
In his landmark  work \cite{Riviere-2007}, Rivi\`ere proposed the
second order linear elliptic  system \begin{equation}\label{eq:Riviere 2007}
	-\Delta u=\Omega\cdot \nabla u \qquad \text{in }B^2\subset\R^2,
\end{equation}
with  $\Omega=(\Omega_{ij})\in L^2(B^2,so_m\otimes \Lambda^1\R^2)$  and $u\in W^{1,2}(B^2, \R^m)$, which models
the Euler-Lagrange equations of  critical points of all second order conformally invariant variational functionals over maps $u\in W^{1,2}(B^2,N)$, where   $B^2\subset \R^2$ is the unit disk, and $N\subset \R^m$ is an arbitrary compact Riemannian manifold. In particular, \eqref{eq:Riviere 2007} includes the equation of weakly harmonic maps from  $B^2$ to $N$ and the prescribed mean curvature equations. A crucial observation of \cite{Riviere-2007} is a conservation law induced by the anti-symmetry of $\Om$, from which the continuity of weak solutions to equation \eqref{eq:Riviere 2007} follows. This gave an affirmative answer to the long standing conjectures of Hildebrandt and Heinz,  and an alternate proof of
Helein's celebrated regularity theorem on weak harmonic maps in dimension two. The technique developed in \cite{Riviere-2007} has also profound applications beyond conformally invariant problems; see \cite{Riviere-2011,Riviere-2012} for a comprehensive overview.

Rivi\`ere and Struwe have further considered in \cite{Riviere-Struve-2008}
the same system as \eqref{eq:Riviere 2007} in supercritical dimensions $n\ge 3$:
\begin{equation}\label{eq:Riviere 2007-2}
	-\Delta u=\Omega\cdot \nabla u \qquad \text{in }B^n\subset\R^n,
\end{equation}
where $\Omega=(\Omega_{ij})\in L^2(B^n,so_m\otimes \Lambda^1\R^n)$.  Although there is no conservation law associated with \eqref{eq:Riviere 2007-2}
for $n\ge 3$, Rivi\`ere and Struwe managed to transform  \eqref{eq:Riviere 2007-2} into a gauge equivalent system through
Uhlenbeck's gauge construction associated with $\Omega$.  It was established in \cite{Riviere-Struve-2008}
that a local H\"older regularity holds for any weak solution $u$ to
\eqref{eq:Riviere 2007-2} under the smallness condition
\[\sup _{x \in B^n_1, r>0}\left(\frac{1}{r^{n-2}} \int_{B^n_{r}(x) \cap B}\left(|\nabla u|^{2}+|\Omega|^{2}\right) d x\right)^{1 / 2}<\varepsilon(n,m).\]
As an application, they reproved the partial regularity theorem on stationary harmonic maps in dimensions $n\ge 3$, due to Evans \cite{Evans-1991} and Bethuel \cite{Bethuel-1993}.

The techniques in \cite{Riviere-2007,Riviere-Struve-2008} have been subsequently extended to fourth order elliptic systems
by Lamm and Rivi\`ere  \cite{Lamm-Riviere-2008} in dimension $n=4$ and Struwe \cite{Struwe-2008} for $n\ge 5$ in the course
of study of biharmonic maps.
Recall that an extrinsic (or intrinsic resp.) biharmonic map from $B^n$ into a closed Riemannian manifold $N$
is a critical point of the energy functional
 \[\int_{B^n}|\De u|^2 \quad \big(\text{or}\ \int_{B^n}|(\De u)^T|^2 \quad \text{ resp.}\big) \qquad  \text{ for } u\in W^{2,2}(B^n, N),\]
 where $(\De u)^T$ is the orthogonal projection of $\De u$ onto the tangent space $T_u N$.
In \cite{Lamm-Riviere-2008}, the authors formulated the following system of 4th order linear elliptic equations
\begin{equation}\label{eq:Lamm-Riviere 2008}
	\De^{2}u=\De(V\cdot\na u)+{\rm div}(w\na u)+F \cdot\na u \quad  \text{in }B^4,
\end{equation}
where $V, w$ belong to certain function spaces and $F=\na \om+W$ 
with  $\om\in L^{2}(B^4,so_m)$ being antisymmetric. By constructing a corresponding conservation law for system \eqref{eq:Lamm-Riviere 2008},
an everywhere continuity for weak solutions  of \eqref{eq:Lamm-Riviere 2008} was established in \cite{Lamm-Riviere-2008}.
The approach of \cite{Lamm-Riviere-2008} was further refined by Guo and Xiang in \cite{Guo-Xiang-2019-Boundary}, where a local H\"older continuity for weak solutions of \eqref{eq:Lamm-Riviere 2008} was proven.
The result of  \cite{Lamm-Riviere-2008}  has been applied  to the  theory of regularity for heat flow of biharmonic maps in dimension four.
In \cite{Struwe-2008}, Struwe revisited biharmonic maps in supercritical dimensions $n\ge 5$ and formulated
the following fourth order linear elliptic system:
\begin{equation}\label{eq: Struwe's equation}
	\De^{2}u=\De(D\cdot\na u)+{\rm div}(E\cdot \na u)+(\De\Om+G)\cdot\na u \qquad \text{in }B^n, \end{equation}
where  $D, E, G$ belong to certain function spaces and  $\Om$ is an $so_m$-valued function with entries in $\R^n$.
We refer interested readers to \cite{Lamm-Riviere-2008,Struwe-2008} for detailed computations of writing the  equation of biharmonic maps in the form of \eqref{eq:Lamm-Riviere 2008} or \eqref{eq: Struwe's equation}.
By extending the approach of Rivi\`ere and Struwe \cite{Riviere-Struve-2008}, Struwe  established in \cite{Struwe-2008}
a partial regularity theory for \eqref{eq: Struwe's equation},  under the  growth condition \eqref{eq:GC fourth order} below, which in turn  gave an alternate proof of the H\"older regularity theorems of Chang, Wang and Yang \cite{Chang-W-Y-1999} and Wang \cite{Wang-2004-CPAM,Wang-2004-MZ} for biharmonic maps.
Because of structural similarities, it seems natural to extend the result of Rivi\`ere and Struwe \cite{Riviere-Struve-2008}
on the system of second order linear equations \eqref{eq:Riviere 2007-2} to system of fourth order linear equations \eqref{eq:Lamm-Riviere 2008} and \eqref{eq: Struwe's equation}. Indeed, Struwe raised the following  question in \cite{Struwe-2008}:

\medskip
\noindent\textbf{Struwe's Question}. \emph{It
would be interesting to see if our method can be extended to general linear systems of fourth order that exhibit a structure similar to the one of equation \eqref{eq: Struwe's equation}, as is the case for second  order systems \eqref{eq:Riviere 2007-2}, or in the ``conformal'' case $n= 4$ considered in \cite{Lamm-Riviere-2008}. }
\medskip

Struwe's Question in  the ``conformal'' case $n=4$ has recently been solved  by Guo and Xiang in \cite{Guo-Xiang-2019-Higher}.
More precisely, it was proven in \cite{Guo-Xiang-2019-Higher} that in critical dimensions $n=2k$ for any $k\ge 2$,
a H\"older continuity holds for any weak solution $u\in W^{k,2}(B^n,\R^m)$ of the $2k$-order linear elliptic system with antisymmetric potentials introduced by de Longueville and Gastel in
\cite{deLongueville-Gastel-2019}.  \cite{Guo-Xiang-2019-Higher} was built upon  the ideas by  Rivi\`ere-Struwe \cite{Riviere-Struve-2008} and utilized both Uhlenbeck's gauge transformation and the duality of Lorentz spaces $L^{p,1}-L^{p^\prime, \wq}$, where $1<p<\wq $ and $p^\prime=p/(p-1)$. However, when dimensions $n\ge 5$,  the approach by \cite{Guo-Xiang-2019-Higher}  (see \cite[Section 5]{Guo-Xiang-2019-Higher})
encountered serious technical difficulties, which left open Struwe's Question
in supercritical dimensions $n\ge 5$. Another interesting problem, closely related to the regularity theory on
\eqref{eq:Riviere 2007-2} and Struwe's Question on \eqref{eq: Struwe's equation},
is to study the corresponding inhomogeneous system
of \eqref{eq: Struwe's equation} in dimensions $n\ge 4$. These problems lead us to ask

\smallskip
\begin{problem}\label{problem:TAMS}
	Establish a $L^p$-regularity theory for weak solutions of the fourth order inhomogeneous  elliptic system of Lamm and Rivi\`ere  \cite{Lamm-Riviere-2008} or  Struwe \cite{Struwe-2008}
    $$\De^{2}u=\De(D\cdot\na u)+{\rm div}(E\na u)+F\cdot\na u+f \qquad \text{in }B^n$$
in dimensions  $n\ge 4$. \end{problem}
\smallskip

More specifically, Problem \ref{problem:TAMS} asks that for $f\in L^p(B^n,\R^m)$ with $1<p<\infty$,
if a $W^{4,p}_{\loc}$-regularity holds for  weak solutions of the linear systems \eqref{eq:Lamm-Riviere 2008} or \eqref{eq: Struwe's equation}, provided  certain smallness conditions are imposed on both the linear coefficient functions and the solution.
In the critical dimension $n=4$, Problem \ref{problem:TAMS} was solved by Guo, Xiang and Zheng in \cite{Guo-Xiang-Zheng-2021-CV},
where they proved that  if $f\in L^p$ for $1<p<4/3$, then $u\in W^{3, 4p/(4-p)}_\loc\subset C^{0, 4(1-1/p)}_\loc$. In particular, this implies that  when $n=4$, every weak solution of the system \eqref{eq:Lamm-Riviere 2008} or \eqref{eq: Struwe's equation} is locally $\al$-H\"older continuous for all $0<\al<1$. A similar $L^p$-theory for general even order linear elliptic systems proposed by de Longueville and Gastel \cite{deLongueville-Gastel-2019} was also established by \cite{Guo-Xiang-Zheng-2021-Lp}
in critical dimensions. For  applications  to  biharmonic maps, see  Laurain-Lin  \cite{Laurain-Lin-2021-Crelle} for energy convexity
and Laurain-Rivi\`ere \cite{Laurain-Riviere-2013-ACV} and Wang-Zheng \cite{Wang-Zheng-2012-JFA} for energy quantization. We also point out that the theory of biharmonic maps has been successfully applied in Cheng-Zhou's solution of the Rosenberg-Smith conjecture in their recent work \cite{Cheng-Zhou-2021}.
We would like to mention that a positive answer to Problem \ref{problem:TAMS} would solve Struwe's Question.
However, Problem \ref{problem:TAMS} remains open in supercritical dimensions $n\ge 5$. In this paper,
we will make some partial progress towards Problem \ref{problem:TAMS}.

In the second order case, motivated by the study on approximate harmonic maps and
heat flow of harmonic maps,  we would like to ask the following  problem.
\medskip

	\emph{Develop a  $W^{2,p}$-regularity theory for  the inhomogeneous Rivi\`ere's system
	\begin{equation}\label{eq:inhomogenuous system ST}
		-\Delta u=\Omega\cdot \nabla u+f\qquad \text{in }B^n,
	\end{equation}
   where $\Omega\in L^2(B^n, so_m\otimes\R^n)$  and $f\in L^p(B^n, \R^m)$.}
\medskip

 This problem was first considered by Sharp and Topping \cite{Sharp-Topping-2013-TAMS} in dimension $n=2$.
 Utilizing the conservation law of Rivi\`ere \cite{Riviere-2007}, they proved that if $f\in L^p(B^2,\R^m)$ for  $p\in (1,2)$, then every weak solution $u\in W^{1,2}(B^2,\R^m)$  belongs to $W^{2,p}_{\loc}(B^2,\R^m)\subset C^{0,2(1-1/p)}_{\loc}$. In particular, any weak solution of \eqref{eq:Riviere 2007} is locally $\alpha$-H\"older continuous for any $0<\alpha<1$.
See Laurain-Rivi\`ere \cite{Laurain-Riviere-2014-APDE} and Lamm-Sharp \cite{Lamm-Sharp-2016-CPDE} for some further related results.

For dimensions $n\ge 3$, in the course of studying the heat flow of harmonic maps,  Moser \cite{Moser-2015-TAMS} considered 
the $L^p$-regularity of the system of  approximate harmonic maps $u:B^n\to N$:
\begin{equation}\label{moser}
-\De u=A(u)(\na u, \na u)+f,
\end{equation}
and proved that, if $f\in L^p(B^n,\R^m)$ for some $n/2<p<\wq$, then $u\in W^{2,p}_\loc(B^n, N)$, under certain  smallness
condition on $\nabla u$.
One crucial idea  of \cite{Moser-2015-TAMS} is to rewrite the system \eqref{moser} via  the Gauge transformation of Rivi\`ere and Struwe \cite{Riviere-Struve-2008}. On the other hand,  Sharp \cite{Sharp-2014} established a Morrey-space regularity
for the linear system \eqref{eq:inhomogenuous system ST}, namely,
$M^{\frac{2p}{n}, n-2}$-regularity for $\nabla^2 u$ holds under a smallness condition
on $\|\Omega\|_{M^{2,n-2}}$.
In view of Moser \cite{Moser-2015-TAMS}, Sharp made the following expectation in \cite[Remark 1.3]{Sharp-2014}:
\smallskip

\noindent\textbf{Sharp's expectation.}
\emph{One would expect Moser's $L^p$-regularity on \eqref{moser}
remains to hold for the system  \eqref{eq:inhomogenuous system ST},
under the additional condition
\begin{equation*}\tag{GC-2}\label{eq:growth condition second order}
	|\Om|\le C|\na u|.
\end{equation*}}
We will give an affirmative answer to  this expectation in Theorem \ref{thm:second order} below.	
	
\subsection{Main results}	

Henceforth, we will assume $m>1, n\geq 5$. Let $B_{r}=\{x\in\R^{n}:|x|<r\}$ and $u\in W^{2,2}(B_{2},\R^{m})$.
Consider the following inhomogeneous 4th order elliptic system
\begin{eqnarray}\label{Struwe_type}
	\De^{2}u=\De(D\cdot\na u)+{\rm div}(E\na u)+F\cdot\na u+f &  & \text{in }B_{2},\label{eq: Struwe system}
\end{eqnarray}
with $F=\De\Om+G$, and
\begin{equation}\label{eq: regularity of coefficients}
\begin{aligned}
 &   D\in W^{1,2}(B_{2},M_{m}\otimes\Lambda^{1}\R^{n}),\quad\quad E\in L^{2}(B_{2},M_{m})\\
 &  \Om\in W^{1,2}(B_{2},so_{m}),\quad\quad   G\in L^{\frac{4}{3},1}(B_{2},M_{m}\otimes\Lambda^{1}\R^{n}).
\end{aligned}
\end{equation}
In coordinates, \eqref{Struwe_type} reads as
 \[\De^2 u^i= \De (D_j^i\cdot \na u^j)+{\rm div}(E^i_j \na u^j)+F_{j}^i\cdot \na u^j+f^i, \quad 1\le i\le m,\]
where the Einstein summation convention is used for repeated indices.

In this paper, we aim to establish an $L^p$-regularity theory for \eqref{eq: Struwe system} under the following growth condition
on $D, E, G, \Om$:
 \begin{equation}
	\begin{aligned}|D|+|\Omega| & \leq C|\nabla u|,\\
		|E|+|\nabla D|+|\nabla\Omega| & \leq C\left|\nabla^{2}u\right|+C|\nabla u|^{2},\\
		|G| & \leq C\left|\nabla^{2}u\|\nabla u\right|+C|\nabla u|^{3}.
	\end{aligned}
	\tag{GC-4}\label{eq:GC fourth order}
\end{equation}

Although the $L^p$-theory of  \eqref{eq: Struwe system} under condition (GC-4)
does not answer Problem 1.1, it provides very interesting insights on attacking this challenging problem.
From the analytic point of view, the nonlinearity under (GC-4) is of critical growth so that for a weak solution $u\in W^{2,2}(B^n,\R^m)$ there merely holds $|\na^2 u|^2\in L^1(B^n)$
and  the standard $L^p$-regularity theory is not applicable. Furthermore, 
the nonlinearity is so strong that it is also impossible to apply the standard bootstrapping argument,
even if some improved regularity, e.g. $\na^2 u\in L^{2+\ep}$ for a sufficiently small $\ep>0$, is assumed.

Since \eqref{Struwe_type} models biharmonic maps when $f\equiv 0$
(see Lamm and Rivi\`ere  \cite{Lamm-Riviere-2008} and Struwe \cite{Struwe-2008}),
our $L^p$-regularity theory, via the Sobolev embedding theorem,  implies that
$u\in C^\alpha_{\rm{loc}}(B^n)$ for any $\alpha\in (0,1)$, which in turn improves
Struwe's H\"older regularity theorem.  Note that the system for both approximate biharmonic maps
and heat flow of biharmonic maps do satisfy both \eqref{Struwe_type} and the growth condition \eqref{eq:GC fourth order},
hence the $L^p$-regularity theory will have direct application to the study of biharmonic maps
and heat flow of biharmonic maps in supercritical dimensions.

We denote by $M^{p,\la}(B_r)$ and $M_\ast^{p,\la}(B_r)$ the $(p,r)$-Morrey space and weak $(p,r)$-Morrey space respectively (see Section \ref{sec: preliminary} for
their definitions). Our first  theorem is stated as follows.
\begin{theorem} \label{thm: Morrey decay estimate}
   Suppose $f\in M^{1,n-4+\alpha}(B_2,\R^m)$ for some $\al\in (0,1)$ and $u\in W^{2,2}(B_{2},\R^{m})$ is a solution of system \eqref{eq: Struwe system} satisfying \eqref{eq:GC fourth order}.
	There exist constants $\ep=\ep(m,n,\al)$ and $C=C(m,n,\al)>0$ such that if
	\begin{equation}
		\|\na^{2}u\|_{M^{2,n-4}(B_{1})}+\|\na u\|_{M^{4,n-4}(B_{1})}\le\ep,\label{eq: smallness assumption}
	\end{equation}
	then
	\begin{eqnarray*}
		\na u\in M_{\ast}^{4,n-4+4\al}(B_{1/2}) & \text{and} & \na^{2}u\in M_{\ast}^{2,n-4+2\al}(B_{1/2}).
	\end{eqnarray*}
	Moreover,
	\begin{equation}
		\|\na u\|_{M_{\ast}^{4,n-4+4\al}(B_{1/2})}\le C\left(\|\na u\|_{M_{\ast}^{4,n-4}(B_{1})}+\|f\|_{M^{1,n-4+\alpha}(B_{1})}\right),\label{eq: 1 order decay estimate}
	\end{equation}
	\[
	\|\na^{2}u\|_{M_{\ast}^{2,n-4+2\al}(B_{1/2})}\le C\left(\|\na^{2}u\|_{M_{\ast}^{2,n-4}(B_{1})}+\|\na u\|_{M_{\ast}^{4,n-4}(B_{1})}+\|f\|_{M^{1,n-4+\alpha}(B_{1})}\right).
	\]
\end{theorem}

Theorem \ref{thm: Morrey decay estimate}, relaxing the $L^p$-assumption on $f$ to a Morrey assumption on $f$, seems to be new even in the critical dimension $n=4$.
In connection with the Morrey smallness assumption \eqref{eq: smallness assumption},
the Morrey assumption on $f$ seems to be more compatible than the  $L^p$ assumption on $f$.

The smallness assumption \eqref{eq: smallness assumption} is natural in terms of both translation and dilation
invariance.  When $f\equiv 0$,
the monotonicity formula for stationary biharmonic maps justifies this smallness assumption, see e.g. Wang \cite{Wang-2004-CPAM}, Struwe \cite{Struwe-2008} and Moser \cite{Moser-2008-CPDE}. For heat flow of biharmonic map flow (i.e. $f=u_t$),
a parabolic version of smallness assumptions can also be verified in some cases, see e.g. Hineman-Huang-Wang \cite{Hineman-Huang-Wang-2014}.
In view of the embedding $M_{\ast}^{4,n-4+4\al}(B_1)\subset M^{q,n-q+q\al}(B_1)$
for any $1\le q<4$, the estimate \eqref{eq: 1 order decay estimate}
can be viewed as a slight improvement of Struwe \cite[Estimate (37)]{Struwe-2008},
where it was proved that  $\na u\in M^{q,n-q+q\al}(B_{1/2})$ for some $1<q<2$.

An immediate consequence of Theorem \ref{thm: Morrey decay estimate} is the following optimal H\"older regularity.
\begin{corollary}\label{coro:Holder}
	 Suppose $f\in M^{1,n-4+\alpha}(B_2,\R^m)$ for some $\al\in (0,1)$ and $u\in W^{2,2}(B_{2},\R^{m})$ is a solution of system \eqref{eq: Struwe system} satisfying \eqref{eq:GC fourth order}.  There exist constants $\ep=\ep(m,n,\al)$ and $C=C(m,n,\al)>0$ such that  if \eqref{eq: smallness assumption} holds,
	then $u\in C^{0,\al}_{\loc}({B}_{1})$ with
	\[
	\|u\|_{C^{0,\al}({B}_{1/2})}\le C\left(\|\na u\|_{M_{\ast}^{4,n-4}(B_{1})}+\|f\|_{M^{1,n-4+\alpha}(B_{1})}\right).
	\]
\end{corollary}

When $f\equiv 0$, this implies that solutions of \eqref{eq: Struwe's equation} are locally $\alpha$-H\"older continuous with any exponent $0<\alpha<1$, which improves the main result of Struwe \cite{Struwe-2008}. See Rupflin \cite{Rupflin-2008} and Wang-Zheng \cite{Wang-Zheng-2012-JFA} for some related results.

Theorem \ref{thm: Morrey decay estimate} provides the key technical tool  to prove  the following $L^p$-regularity result.

\begin{theorem}\label{thm: main result}
Suppose $f\in L^{p}(B_{1})$ for some $n/4<p<\wq$ and $u\in W^{2,2}(B_{1},\R^{m})$ is a solution of system \eqref{eq: Struwe system} satisfying  \eqref{eq:GC fourth order}.
\begin{itemize}
	\item[(i)]  When $p<n$, there exists a constant $\ep=\ep(m,n,p)$ such that
	if the assumption \eqref{eq: smallness assumption} holds,
	then $u\in W_{\loc}^{3,\frac{np}{n-p}}(B_{1})$ and \[
\|\na^{3}u\|_{L^{\frac{np}{n-p}}(B_{1/2})} \le C_*\left(\|\na u\|_{M^{4,n-4}(B_{1})}+\|\na^{2}u\|_{M^{2,n-4}(B_{1})}
+\|f\|_{L^{p}(B_{1})}\right).
\]
Here $C_*=c_*(1+\ep+\|f\|_{L^p(B_1)})^{a_*}$ for constants $c_*$ and $a_*$ depending on $n, m, p$.
	
	\item[(ii)] When $p\ge n$, for any $1<q<\wq$, there exists a constant $\ep=\ep(m,n,p,q)$  such that
	if the smallness assumption \eqref{eq: smallness assumption} holds, then $	u\in	W_{\loc}^{3,q}(B_{1})$ and \[
\|\na^{3}u\|_{L^{q}(B_{1/2})} \le D_*\left(\|\na u\|_{M^{4,n-4}(B_{1})}+\|\na^{2}u\|_{M^{2,n-4}(B_{1})}
+\|f\|_{L^{p}(B_{1})}\right).
\]
Here $D_*=d_*(1+\ep+\|f\|_{L^p(B_1)})^{b_*}$ for constants $d_*$ and $b_*$ depending on $n, m, p, q$.

\end{itemize}
\end{theorem}

We would like to remark that in Theorem \ref{thm: main result}, the constant $C_*$
(and $D_*$) depends not only on $\ep, n,m,p$ (and $q$) but also on $\|f\|_{L^p(B_1)}$.
The reason for dependence on $\|f\|_{L^p(B_1)}$ is that our system is nonlinear and in general one can only expect a priori estimates involving polynomial dependency on $\|\na u\|_{M^{4,n-4}(B_{1})}$, $\|\na^{2}u\|_{M^{2,n-4}(B_{1})}$ and $\|f\|_{L^{p}(B_{1})}$, under the Morrey smallness assumption \eqref{eq: smallness assumption} on $\nabla u$ and $\nabla^2 u$.

This  result provides an affirmative answer to Problem \ref{problem:TAMS} under the growth condition \eqref{eq:GC fourth order}. As a simple consequence of this theorem and the Sobolev embedding theorem, we can infer that the smallness assumption \eqref{eq: smallness assumption} implies
\[
u\in\begin{cases}
	C_{\loc}^{0,\al}(B_{1}), & \text{if }n/4<p<n/3,\\
	C_{\loc}^{1,\al-1}(B_{1}), & \text{if }n/3<p<n/2,\\
	C_{\loc}^{2,\al-2}(B_{1}), & \text{if }n/2<p<n,
\end{cases}
\]
where $\al=4-n/p$,  whenever $p<n$.  As a further application, we have

\begin{corollary}\label{thm: main result 2} Suppose $f\in L^{p}(B_{1})$ for some $ n/4<p<\wq$ and $u\in W^{2,2}(B_{1},\R^{m})$ a solution of system \eqref{eq: Struwe system} satisfying  \eqref{eq:GC fourth order}. Suppose further that  $$
|\nabla D|\leq C(|\nabla^2 u|+|\nabla u|^2),$$
$$|\nabla E|+|\nabla^2 \Omega|\leq C(
|\nabla^3 u|+|\nabla^2 u||\nabla u|+|\nabla u|^3).
$$
Then there exists a constant $\ep=\ep(m,n,p)$ such that
if the assumption \eqref{eq: smallness assumption} holds, then $u\in W_\loc^{4,p}(B_1)$ and
\[
\|\na^{4}u\|_{L^{p}(B_{1/2})} \le C\left(\|\na u\|_{M^{4,n-4}(B_{1})}+\|\na^{2}u\|_{M^{2,n-4}(B_{1})}
+\|f\|_{L^{p}(B_{1})}\right).
\]
Here $C=c(1+\ep+\|f\|_{L^p(B_1)})^{a}$ for two constants $c$ and $a$ depending on $n, m, p$.
In particular, any  approximate  biharmonic map $u$
with  drift term $f\in L^{p}(B_{1})$ for some $ n/4<p<\wq$  belongs to  $W^{4,p}_\loc (B_1),$ provided the smallness condition  \eqref{eq: smallness assumption} holds for a sufficiently small $\ep$.
\end{corollary}

Corollary \ref{thm: main result 2} extends the corresponding results  by Wang-Zheng \cite{Wang-Zheng-2012-JFA} and Laurain-Rivi\`ere  \cite{Laurain-Riviere-2013-ACV} in the critical dimension $n=4$. It also leads to the following weak compactness and energy gap results.
\begin{corollary}\label{thm: Weak compactness}
	There is a sufficient small constant $\ep=\ep(m,n)>0$ such that
\begin{description}
	\item[(1) Weak Compactness] For any sequence  $u_k\in W^{2,2}(B_{1},N)$  of biharmonic maps  which converges weakly to a map
	$ u\in W^{2,2}(B_{1},N)$,  if \[	\|\na^{2}u_k\|_{M^{2,n-4}(B_{1})}+\|\na u_k\|_{M^{4,n-4}(B_{1})}\le\ep, \qquad \forall\,\,k\ge 1,	\] then up to a subsequence,  $u_k\to u$ strongly in $W^{2,2}_{\rm{loc}}(B_1, N)$. In particular, $u$ is a smooth biharmonic map; and
	
	\item[(2) Energy gap] If $u\in W^{2,2}(\R^n, N)$ is a  biharmonic map  satisfying \[	\|\na^{2}u\|_{M^{2,n-4}(\R^n)}+\|\na u\|_{M^{4,n-4}(\R^n)}\le\ep, \]
	then $u\equiv p$ in $\R^n$ for a point $p\in N$.
\end{description}
\end{corollary}
	
For geometric  applications of this type of result, see Wang-Zheng \cite{Wang-Zheng-2012-JFA} and Laurain-Rivi\`ere  \cite{Laurain-Riviere-2013-ACV}
on the energy identity of biharmonic maps in  dimension $n=4$.

The requirement $p>{n}/{4}$ in Theorem \ref{thm: main result} ensures that $f\in M^{1,n-4+\al}(B_1)$ for some $0<\al<1$. It is natural to ask what happens when $1<p\le {n}/{4}$. Observe that the heat flow $u$ of biharmonic maps  can be viewed as \eqref{eq: Struwe system} for $f=u_t\in L^p$, with $p=2\le n/4$ when dimensions $n\ge 8$, see e.g.  Moser \cite{Moser-2009-ACV}.
This motivates us to consider the case $f\in L^p$ for $1<p\le n/4$. We can prove
\begin{theorem}\label{thm:p small}
Suppose $f\in L^{p}\cap M^{1,n-4+\al}(B_1)$ for
some $1<p\le n/4$ and $0<\al<1$ and $u\in W^{2,2}(B_{1},\R^{m})$
is a weak solution to system \eqref{eq: Struwe system} satisfying  \eqref{eq:GC fourth order}. There
exists $\ep>0$ such that if the smallness condition \eqref{eq: smallness assumption}
holds, then
\[
\na^{2}u\in L^{p\eta}\cap M_{\ast}^{\eta,n-\eta(2-\al)}(B_{1/2}) \quad \text{and}
\quad
\na u\in L^{p\eta\chi}\cap M_{\ast}^{\eta\chi,n-\eta\chi(1-\al)}(B_{1/2}),
\]
where $\chi=(2-\al)/(1-\al)>2$ and $\eta=(4-\al)/(2-\al)>2$.
\end{theorem}

We would like to remark that with slight changes of arguments, all results stated as above remain to hold
if the coefficient function $F$ in equation \eqref{eq: Struwe system}  takes
the form $F=\na \om+W$ of equation \eqref{eq:Lamm-Riviere 2008} and (GC-4) is replaced by a corresponding one.

Finally, as aforementioned, as a by-product of our method, we provide an affirmative answer to Sharp's expectation.
\begin{theorem}\label{thm:second order}
	Suppose $f\in L^p(B_1)$ for some   $1\le {n}/{2}<p<\infty$ and $u\in W^{1,2}(B_1,\R^m)$ is a weak solution of system \eqref{eq:inhomogenuous system ST}. If, in addition, $\Om\in L^{2}(B_1,so_{m}\otimes\wedge^{1}\R^{n})$ satisfies the growth condition \eqref{eq:growth condition second order},
	 then there exists $\epsilon=\epsilon(m,n,p)>0$ such that
	$$u\in W^{2,p}_{\loc}(B_1,\R^m),$$
	whenever  $\|\nabla u\|_{M^{2,n-2}(B_1,\R^m)}<\epsilon$.
	
\end{theorem}
	
Theorem \ref{thm:second order}  extends the main result of Moser \cite[Theorem 1.2]{Moser-2015-TAMS} for $p>n/2$ with an
alternate proof. It is an interesting question that if Theorem  \ref{thm:second order} holds
when $1<p\le n/2$.

\subsection{Strategy of the proof}	

To derive Theorem \ref{thm: Morrey decay estimate},  we shall first rewrite the system using the Gauge transform  of  Struwe \cite{Struwe-2008}, and then apply the Hodge decomposition to simplify the problem. Morrey type decay estimates then follow from a combination of Riesz potential theory  and a decay property of harmonic functions.

 Theorem \ref{thm: main result} follows from a delicate iteration argument. To explain this strategy clearly,
 we first sketch the proof of Theorem \ref{thm:second order}.

\begin{proof}[Proof of Theorem \ref{thm:second order}]
We first consider the case $\frac{n}{2}<p<n$. Sharp \cite[Theorem 1.2]{Sharp-2014} has proved that $\na u\in M_{\loc}^{2,n-2+2\al}(B_1)$ with $\alpha=2-{n}/{p}\in (0,1)$. It follows from the
growth condition that $\Om\in M_{\loc}^{2,n-2+2\al}(B_1)$. Therefore
$\Om\cdot\na u\in M_{\loc}^{1,n-2+2\al}(B_1)$.

First suppose $\al<1/2$. Extend $\Om$, $u$ and $f$ from $B_{1/2}$ into $\R^{n}$ with compact
support in $B_{2}^{n}$ in a norm-bounded way. Let $u_{1}=I_{2}(\Om\cdot\na u)$
and $u_{2}=I_{2}(f)$ such that $h=u-u_{1}-u_{2}$ is a harmonic function
in $B_{1/2}$, where  $I_\al=c|x|^{\al-n}$ is the standard Riesz potential. Since $\Om\cdot\na u\in M^{1,n-2+2\al}\cap L^{1}(\R^{n})$,
Adams' potential theory (see Proposition \ref{prop: Riesz in Morrey} below) implies
\[
|\na u_{1}|\le CI_{1}(\Om\cdot\na u)\in L^{2\chi,\wq}(\R^{n}),
\]
where
\[
\chi=\frac{1}{2}\left(\frac{2-2\al}{1-2\al}\right)>1.
\]
By standard elliptic regularity theory, $u_{2}\in W^{2,p}(\R^{n})$. Hence
$\na u_{2}\in L^{p^*}(\R^{n})$, where $p^*=\frac{np}{n-p}$.

Since $h\in C^\wq$, if $2\chi>p^{\ast}$,
then we find that $\na u\in L_{\loc}^{p^{*}}(B_1)$.
If $2\chi\le p^{\ast}$, we obtain $\na u\in L^{2\chi,\wq}(B_{1/2})$.
Since $2\chi>2$, we have $\na u\in L_{\loc}^{q}(B_{1})$ for some
$q>2$. Then the growth condition \eqref{eq:growth condition second order} implies that $\Om\in M_{\loc}^{2,n-2+2\al}\cap L^{q}_\loc (B_1)$.
It follows that
\[
\Om\cdot\na u\in M_{\loc}^{1,n-2+2\al}(B^{n})\cap L_{\loc}^{q/2}.
\]
Since $q/2>1$, using Adams' potential (see Proposition \ref{prop: Adams-Riesz}) again gives
$
\na u_{1}\in L^{\chi q}.
$
If $\chi q\le p^{\ast}$, we obtain
$
\na u\in L^{\chi q}.
$
Thus, we find the iteration:
\[
\na u\in L^{q}\Rightarrow\na u\in L^{\chi q}.
\]
Since $\chi>1$, we can assume that $\chi^{k}q\le p^{\ast}<\chi^{k+1}q$
for some $k\ge1$. After finitely many times iteration, we find that
$\na u\in L_{\loc}^{p^{*}}(B_{1})$.

In the case $\al\ge1/2$, we use embedding $M^{2,n-2+2\al}(B_{1})\subset M^{2,n-2+2\be}(B_{1})$
for any $\be<1/2$ so as to obtain the same regularity as in the case
$\al<1/2$. As a consequence, we can always  derive $\na u\in L_{\loc}^{p^{*}}(B_{1})$. Now the second order regularity $u\in W^{2,p}_{\loc}(B_1)$ follows from the usual elliptic regularity theory.

If $p\geq n$, then $f\in L^q(B^n)$ for any ${n}/{2}<q<n$ . Running the previous argument we conclude that $\nabla u\in L^{\frac{nq}{n-q}}_{\loc}(B_1)$. This implies that $u_1\in \bigcap_{1<q<\wq} W^{1,q}_{\loc}(B_1)$ and so  finally $u\in W^{2,p}_{\loc}(B_1)$. The proof is complete.
\end{proof}

Our proof of Theorem \ref{thm: main result} follows a similar approach, but the analysis becomes much more involved. In a first step, we derive the weak Morrey decay estimate for solutions of system \eqref{eq: Struwe system}, that is, Theorem \ref{thm: Morrey decay estimate}.  Unlike the case of linear systems in \cite{Guo-Xiang-Zheng-2021-CV,Guo-Xiang-Zheng-2021-Lp,Sharp-2014,Sharp-Topping-2013-TAMS}, this regularity improvement is not strong enough for iteration yet. To fill the gap,  two observations are needed here:
\begin{itemize}
	\item The  weak Morrey regularity of $\na^2 u$ automatically implies an improvement of $\na u$, i.e,
	$\na u\in L^{2\chi}\cap M_{\ast}^{2\chi,n-2\chi(1-\al)}(B_{1/4})$, where
	$\chi\equiv(2-\al)/(1-\al)>2$;
	
	\item The growth condition implies a corresponding regularity improvement for the  Gauge transformation (see Lemma \ref{lem: Gauge transform-general-version} below).
\end{itemize}
By the first observation, we obtained an improved regularity of $\na u$.  But this improvement itself is not sufficient for the iteration method yet. To proceed, our new idea is to further track and  improve the regularity of  Gauge transforms. We then turn to
 construct an associated Gauge transform on  smaller balls (half radius of the previous one) with improved regularity. This is realized by the second observation.  Then, combining these improved Gauge transforms, and tracking both the Lebesgue integrability and Morrey regularity of $\na u$ and $\na^2 u$ simultaneously,   an application of the Riesz potential theory gives a further improvement on integrability of $\nabla u$ and $\nabla^2 u$.  Finally, to obtain the optimal interior regularity, we run an iteration scheme by  repeatedly constructing the gauge transforms on a sequence of shrinking balls and then using the gauge equivalent equations on shrinking balls. Surprisingly,  in contrast to those infinite iteration  on linear systems in \cite{Guo-Xiang-Zheng-2021-CV,Guo-Xiang-Zheng-2021-Lp,Sharp-2014,Sharp-Topping-2013-TAMS}, our  iteration process actually only takes finitely many steps thanks to the nonlinearity of the problems.  Finally we mention that a crucial harmonic analysis theory  used in the proof is the boundedness of Riesz operators between weak Morrey spaces,
which is due to Ho \cite{Ho-2014} (see also Proposition \ref{prop: Ho} below).

 Our notations are standard. By $A\lesssim B$, we mean there is an absolute constant $C>0$ such that $A\le CB$. The constant $C$ may differ from line to line.

\section{Preliminaries}\label{sec: preliminary}
In this section, we introduce some function spaces and the related Riesz
potential theory between these function spaces. They play a central
role in later proofs.

\subsection{Morrey spaces}

Let $\Om\subset\R^{n}$ be a smooth domain. For $1\le p<\wq$, let
$L^{p}(\Om)$ be the usual $L^p$ space on $\Om$ and $L_{\ast}^{p}(\Om)$
the weak $L^{p}$ space on $\Om$.

Let $1\le p<\wq$ and $0\le s\le n$. The Morrey space $M^{p,s}(\Om)$
consists of functions $f\in L^{p}(\Om)$ such that
\[
\|f\|_{M^{p,s}(\Om)}\equiv\sup_{x\in\Om,0<r<\diam(\Om)}r^{-s/p}\|f\|_{L^{p}(B_{r}(x)\cap\Om)}<\wq.
\]
The weak Morrey space $M_{\ast}^{p,s}(\Om)$ consists of functions
$f\in L_{\ast}^{p}(\Om)$ such that
\[
\|f\|_{M_{\ast}^{p,s}(\Om)}\equiv\sup_{x\in\Om,0<r<\diam(\Om)}r^{-s/p}\|f\|_{L_{\ast}^{p}(B_{r}(x)\cap\Om)}<\wq.
\]

Note that $M^{p,0}(\Om)=L^{p}(\Om)$ and $M^{p,n}(\Om)=L^{\wq}(\Om)$,
and $M_{\ast}^{p,0}(\Om)=L_{\ast}^{p}(\Om)$. When $\Om$ is a bounded
domain, it follows from H\"older's inequality and the simple embedding
$L_{\ast}^{p}(\Om)\subset L^{q}(\Om)$ ($1\le q<p$) that,

\[
L^{p}(\Om)\subset M^{q,n(1-\frac{q}{p})}(\Om),\qquad\forall\;1\le q<p
\]
and
\[
M_{\ast}^{p,s}(\Om)\subset M^{1,n+\frac{s-n}{p}}(\Om),\qquad\forall\;1<p<\wq.
\]

We shall need the following well-known H\"older's inequality for weak $L^{p}$ functions.

\begin{proposition}\label{prop: Lorentz-Holder inequality}
Let $1<p_{1},p_{2}<\wq$ be such that $\frac{1}{p}=\frac{1}{p_{1}}+\frac{1}{p_{2}}\le1$.
Then, $f\in L_{\ast}^{p_{1}}(\Om)$ and $g\in L_{\ast}^{p_{2}}(\Om)$
implies $fg\in L_{\ast}^{p}(\Om)$. Moreover,
\[
\|fg\|_{L_{\ast}^{p}(\Om)}\le\|f\|_{L_{\ast}^{p_{1}}(\Om)}\|g\|_{L_{\ast}^{p_{2}}(\Om)}.
\]
\end{proposition}

The following proposition concerns H\"older's inequalities in Morrey functions.
The proof is straightforward and thus omitted.

\begin{proposition}\label{prop: Holder for Morrey space} Let $1\le p_{1},p_{2}\le\wq$
and $0\leq q_{1},q_{2}\le n$ be such that
\begin{eqnarray*}
\frac{1}{p}=\frac{1}{p_{1}}+\frac{1}{p_{2}}\le1 & \text{and} & q=\frac{p}{p_{1}}q_{1}+\frac{p}{p_{2}}q_{2}.
\end{eqnarray*}
Then, there hold
\begin{equation}
\|fg\|_{M^{p,q}(\Om)}\le\|f\|_{M^{p_{1},q_{1}}(\Om)}\|g\|_{M^{p_{2},q_{2}}(\Om)}.\label{eq: Morrey-Holder}
\end{equation}
and
\begin{equation}
\|fg\|_{M_{\ast}^{p,q}(\Om)}\le\|f\|_{M_{\ast}^{p_{1},q_{1}}(\Om)}\|g\|_{M_{\ast}^{p_{2},q_{2}}(\Om)}.\label{eq: Morrey-Holder-2}
\end{equation}
\end{proposition}

As we are concerned with  H\"older regularity theory, we need the following weak type
of Morrey's Dirichlet growth theorem.

\begin{proposition}\label{prop: Morrey's DGT} Suppose $\Om$ is
a bounded smooth domain and $u\in L_{\loc}^{1}(\Om)$ such that $\na u\in M_{\ast}^{p,n-p+p\al}(\Om)$
holds for some $1<p<\wq$ and $\al\in(0,1)$. Then $u\in C^{0,\al}(\overline{\Om})$
with
\[
\|u\|_{C^{0,\al}(\overline{\Om})}\le C\|\na u\|_{M_{\ast}^{p,n-p+p\al}(\Om)}
\]
for some $C=C(n,p,\Om)$. \end{proposition}
\begin{proof}
By Poincar\'e's inequality, for any $x\in\Om$ and $0<r<\diam(\Om)$,
there holds
\[
\fint_{B_{r}(x)\cap\Om}|u-u_{r,x}|\le Cr\fint_{B_{r}(x)\cap\Om}|\na u|.
\]
Since $p>1$, we have
\[
\|\na u\|_{L^{1}(B_{r}(x)\cap\Om)}\le Cr^{n(1-1/p)}\|\na u\|_{L_{\ast}^{p}(B_{r}(x)\cap\Om)}\le C\|\na u\|_{M_{\ast}^{p,n-p+p\al}(\Om)}r^{n-1+\al}.
\]
Thus, for any $x\in\Om$ and $0<r<\diam(\Om)$,
\[
\fint_{B_{r}(x)\cap\Om}|u-u_{r,x}|\le C\|\na u\|_{M_{\ast}^{p,n-p+p\al}(\Om)}r^{\al}.
\]
This yields the conclusion by applying Campanato function space theory, see Giaquinta \cite[Chapter III, Theorem 1.2]{Giaquinta-Book}.
\end{proof}
Higher order (weak) Morrey spaces will be useful in our later proofs. For any $k\in\N$,
the $k$th order Morrey space $M_{k}^{p,n-kp}(\Om)$ consists of $f\in W^{k,p}(\Om)$
such that $\na^{l}f\in M^{p,n-lp}(\Om)$ for all $0\le l\le k$,
and we can similarly define the $k$th order weak Morrey space $M_{k,\ast}^{p,n-kp}(\Om)$. It follows from \cite[Proposition 3.2]{Struwe-2008} that $M_{2}^{p,n-2p}(B_1)\subset M_{1}^{2p,n-2p}(B_1)$
with $1<p<n/2$, and
\begin{equation}
\|\na u\|_{M^{2p,n-2p}(B_1)}^{2}\le C\|\na u\|_{M^{1,n-1}(B_1)}\left(\|\na^{2}u\|_{M^{p,n-2p}(B_1)}+\|\na u\|_{M^{p,n-p}(B_1)}\right).\label{Ineq: Struwe}
\end{equation}
In particular, $u\in M_{2}^{2,n-4}(B_1)$ implies that $\na u\in M^{4,n-4}(B_1)$.
Recall that the basic assumption of Struwe \cite{Struwe-2008}
is
\[
R^{n-4}\int_{B_{R}}(|\na^{2}u|^{2}+|\na u|^{4})<\ep,
\]
which together with the monotonicity formula implies that $u\in M_{2}^{2,n-4}(B_{R/2})$
and
\[
\|\na^{2}u\|_{M^{2,n-4}(B_{R/2})}+\|\na u\|_{M^{4,n-4}(B_{R/2})}<C\ep.
\]
Thus, by \eqref{Ineq: Struwe}, one may naturally assume that $u\in M_{2}^{2,n-4}(B_{2})$
satisfies
\[
\|\na^{2}u\|_{M^{2,n-4}(B_{2})}+\|\na u\|_{M^{2,n-2}(B_{2})}<\ep.
\]


We shall frequently use (a special case of) the following Morrey-Sobolev extension\footnote{We would like to thank Prof. Pekka Koskela for pointing out the relevant literatures in this respect.} result due to Burenkov \cite{Burenkov-1975}; see also \cite[Theorem 2.5]{Lamberti-2017} for a new proof.
\begin{proposition}\label{prop:Morrey-Sobolev extension}
For any $k\in \mathbb{N}$, $1\leq p$ and $0\leq s\leq n$, there exists a bounded linear operator $E\colon M^{p,s}_k(B_1)\to M^{p,s}_k(\R^n)$ such that if $f\in M^{p,s}_k(B_1)$, then $Ef=f$ a.e. in $B_1$ and there exists a constant $C=C(k,p,s)>0$ such that for all $f\in M^{p,s}_k(B_1)$, we have
$$\|Ef\|_{M^{p,s}_k(\R^n)}\leq C\|f\|_{M^{p,s}_k(B_1)}.$$
Furthermore, for each $0\leq l\leq k$, there exists a constant $C=C(l,p,s)>0$ such that
$$\|\nabla^l Ef\|_{M^{p,s}(\R^n)}\leq C\|\nabla^l f\|_{M^{p,s}(B_1)}.$$
Similar extension results hold for the higher order weak Morrey-Sobolev spaces $ M^{p,s}_{k,*}(B_1)$ as well.
\end{proposition}
We also refer the interested readers to \cite{Lamberti-2019} for a different construction of the extension operator. Note that in \cite{Lamberti-2019}, the authors only considered the higher order Morrey-Sobolev spaces $M^{p,s}_k(B_1)$. However, the proof works with minor changes (replacing the $L^p$ estimates by corresponding weak $L^p_*$ estimates) for the higher order weak Morrey-Sobolev spaces $M^{p,s}_{k,*}(B_1)$.

\subsection{Riesz potentials}	
	
Let $I_{\al}(x)=c_{\al,n}|x|^{\al-n}$, $0<\al<n$, be the standard
Riesz potentials in $\R^{n}$. The following two propositions are well-known; see Theorem 3.1, Proposition 3.2 and Proposition 3.1 of Adams \cite{Adams-1975}.

\begin{proposition}\label{prop: Riesz in Morrey}  Let $0<\al<n$
and $0\le\la<n$. For $1\le p<(n-\la)/\al$, set
\[
\frac{1}{\tilde{p}}=\frac{1}{p}-\frac{\al}{n-\la}.
\]
Then
\begin{itemize}
	\item[(1)] For every $1<p<(n-\la)/\al$,
	\[
	I_{\al}\colon M^{p,\la}(\R^{n})\to M^{\tilde{p},\la}(\R^{n})
	\]
	is a bounded linear operator;
	
	\item[(2)] For $p=1$,
	\[
	I_{\al}\colon M^{1,\la}(\R^{n})\to M_{\ast}^{\frac{n-\la}{n-\la-\al},\la}(\R^{n})
	\]
	is also a bounded linear operator.
\end{itemize}
\end{proposition}

\begin{proposition} \label{prop: Adams-Riesz} Let $0<\al<\be\le n$
and $1<p<\wq$. Then there exists a constant $C=C_{\alpha,\beta,n,p}>0$
such that for $f\in M^{1,n-\beta}\left(\mathbb{R}^{n}\right)\cap L^{p}\left(\mathbb{R}^{n}\right)$,
there holds
\[
\left\Vert I_{\alpha}f\right\Vert _{\frac{p\beta}{\beta-\alpha},\mathbb{R}^{n}}\leq C\|f\|_{M^{1,n-\beta}\left(\mathbb{R}^{n}\right)}^{\frac{\alpha}{\beta}}\|f\|_{p,\mathbb{R}^{n}}^{1-\frac{\alpha}{\beta}}.
\]
\end{proposition}

In view of the embedding $M_{\ast}^{q,n-q\be}(\R^{n})\subset M^{1,n-\be}(\R^{n})$
for $n/\be\ge q>1$, there holds
\[
\left\Vert I_{\alpha}f\right\Vert _{\frac{p\beta}{\beta-\alpha},\mathbb{R}^{n}}\leq C\|f\|_{M_{\ast}^{q,n-q\beta}\left(\mathbb{R}^{n}\right)}^{\frac{\alpha}{\beta}}\|f\|_{p,\mathbb{R}^{n}}^{\frac{\beta-\alpha}{\beta}}.
\]

Concerning weak Morrey spaces, we will need the following proposition, which is a special case of Ho \cite[Theorem 5.1]{Ho-2014}.

\begin{proposition}\label{prop: Ho} Let $\ensuremath{0<\alpha,\la<n}$
and $1<p<(n-\la)/\al$. Set
\[
\frac{1}{\tilde{p}}=\frac{1}{p}-\frac{\al}{n-\la}.
\]
Then
\[
I_{\al}:M_{\ast}^{p,\la}(\R^{n})\to M_{\ast}^{\tilde{p},\la}(\R^{n})
\]
is a bounded linear operator.
\end{proposition}

As a corollary of Propositions \ref{prop: Adams-Riesz} and \ref{prop: Ho},
for any $\wq>p>1$ and $0<\al<\be<n/p$, we have the following boundedness
result:
\begin{equation}
I_{\al}:M_{\ast}^{p,n-p\be}\cap L^{p}(\R^{n})\to M_{\ast}^{\tilde{p},n-p\be}\cap L^{\tilde{p}}(\R^{n})\qquad\text{where }\tilde{p}=\frac{\be p}{\be-\al},\label{eq: Ho+Adams}
\end{equation}
and
\[
\left\Vert I_{\al}(f)\right\Vert _{L^{\tilde{p}}(\R^{n})}+\left\Vert I_{\al}(f)\right\Vert _{M_{\ast}^{\tilde{p},n-p\be}(\R^{n})}\le C\left(\left\Vert f\right\Vert _{L^{p}(\R^{n})}+\left\Vert f\right\Vert _{M_{\ast}^{p,n-p\be}(\R^{n})}\right).
\]

When the operator under consideration is a singular integral operator,
there holds

\begin{proposition}[Theorem 8.1, \cite{Adams-book-1}]\label{prop: Adams-singular integral operator}
Let $1<p<\wq$ and $0<\la<n$. The usual Calderon-Zygmund singular integral
operators are  bounded on $M^{p,\la}(\R^{n})$.
 \end{proposition}



\section{Morrey estimate  and H\"older continuity}

This section is devoted to prove Theorem \ref{thm: Morrey decay estimate}. For simplicity, denote by $B_{r}=B_{r}(0)\subset\R^{n}$ the open ball centered at origin with radius $r$.
We shall need the following Gauge transform of Struwe \cite[Lemma 3.3]{Struwe-2008}; see also Lamm and Rivi\`ere \cite[Theorem A.5]{Lamm-Riviere-2008} for an equivalent form.

\begin{lemma}[Lemma 3.3, \cite{Struwe-2008}]\label{lem: Gauge transform}
There exist $\ep=\ep(n,m)>0$ and $C=C(n,m)>0$ with the following
property: For every $\Om\in M_{1}^{2,n-4}\cap M^{4,n-4}(B_{1},so_{m}\otimes\wedge^{1}\R^{n})$
with
\[
\|\na\Om\|_{M^{2,n-4}(B_1)}+\|\Om\|_{M^{4,n-4}(B_1)}\le\ep,
\]
there exist $P\in M_{2}^{2,n-4}(B_{1},SO_{m})$ and $\xi\in M_{2}^{2,n-4}(B_{1},so_{m}\otimes\wedge^{n-2}\R^{n})$
such that
\begin{eqnarray}\label{eq: gauge transform equation}
P\D P^{-1}+P\Om P^{-1}=\ast\D\xi &  & \text{in }B_{1},
\end{eqnarray}
and
\begin{eqnarray*}
\D\ast\xi=0\quad\text{in }B_1, &  & \xi=0\quad\text{on }\pa B_{1}.
\end{eqnarray*}
Moreover,
\[
\|\na P\|_{M^{4,n-4}(B_{1})}+\|\na\xi\|_{M^{4,n-4}(B_{1})}\le C\|\Om\|_{M^{4,n-4}(B_{1})}\le C\ep,
\]
\[
\|\na^{2}P\|_{M^{2,n-4}(B_{1})}+\|\na^{2}\xi\|_{M^{2,n-4}(B_{1})}\le C\left(\|\na\Om\|_{M^{2,n-4}(B_{1})}+\|\Om\|_{M^{4,n-4}(B_{1})}\right)\le C\ep.
\]
\end{lemma}

The last two estimates on $P,\xi$ are not separated in the original
statement of Struwe \cite[Lemma 3.3]{Struwe-2008}, but they follow from the proofs there. Below let $P,\xi$ be defined
as in Lemma \ref{lem: Gauge transform}. It follows from the growth
condition \eqref{eq:GC fourth order} on $\Om$ and \eqref{eq: smallness assumption} that
\begin{equation}
\begin{aligned} & \|\na P\|_{M^{4,n-4}(B_{1})}+\|\na\xi\|_{M^{4,n-4}(B_{1})}\le C\|\na u\|_{M^{4,n-4}(B_{1})}\le C\ep.\\
 & \|\na^{2}P\|_{M^{2,n-4}(B_{1})}+\|\na^{2}\xi\|_{M^{2,n-4}(B_{1})}\\
 &\qquad\qquad\qquad\le C\left(\|\na^{2}u\|_{M^{2,n-4}(B_{1})}+\|\na u\|_{M^{4,n-4}(B_{1})}\right)\le C\ep.
\end{aligned}
\label{eq: estimate of gauge transform}
\end{equation}

By \cite[Formula (35)]{Struwe-2008}, the equation of $P\De u$
on $B_{1}$ is given by
\begin{equation}
\De(P\De u)={\rm div}^{2}(D_{P}\otimes\na u)+{\rm div}(E_{P}\cdot\na u)+G_{P}\cdot\na u+\ast\D\De\xi\cdot Pdu+Pf,\label{eq: (35)}
\end{equation}
where  the coefficient functions satisfy the growth condition
\begin{equation}
\begin{aligned}\left|D_{P}\right| & \leq C(|\nabla u|+|\nabla P|),\\
\left|\nabla D_{P}\right|+\left|E_{P}\right| & \leq C\left(\left|\nabla^{2}u\right|+|\nabla u|^{2}+\left|\nabla^{2}P\right|+|\nabla P|^{2}\right),\\
\left|G_{P}\right| & \leq C\left(\left|\nabla^{2}u\right|+\left|\nabla^{2}P\right|\right)(|\nabla u|+|\nabla P|)+C\left(|\nabla u|^{3}+|\nabla P|^{3}\right).
\end{aligned}
\label{eq: growth conditions 2}
\end{equation} For details, see the formula (36)
of \cite{Struwe-2008}.

\begin{proof}[Proof of Theorem \ref{thm: Morrey decay estimate}]\textbf{
}First apply the Hodge decomposition to derive
\begin{eqnarray*}
Pdu=d\tilde{u}_{1}+d^{\ast}\tilde{u}_{2}+\tilde{h} &  & \text{in }B_{1},
\end{eqnarray*}
where $d^{\ast}\tilde{u}_{1}=0$, $d\tilde{u}_{2}=0$ and $\tilde{h}$
is a harmonic 1-form. Note that $\De^{2}\tilde{u}_{1}=\De d^{\ast}(Pdu)$,
$-\De\tilde{u}_{2}=dP\wedge du$ and $\De\tilde{h}=0$ on $B_{1}$.

Next, we extend all the related functions $u$, $\xi$, $P$ and $D_{P}$,
$E_{P}$ and $G_{P}$ from $B_{1}$ into the whole space $\R^{n}$
with compact supports in $B_{2}$ in the same function space in a bounded
way. Set $f\equiv0$ on $B_{1}^{c}$. For simplicity, we keep using the same notations for the extended functions. Then we  define
\begin{equation}
u_{11}=I_{4}\Big({\rm div}^{2}(D_{P}\otimes\na u)+{\rm div}(E_{P}\cdot\na u)+G_{P}\cdot\na u+\ast\D\De\xi\cdot Pdu+\De(\na P\na u)\Big),\label{eq: definition of u11}
\end{equation}
\begin{equation}
u_{12}=I_{4}(Pf),\label{eq: def of u12}
\end{equation}
where $I_{4}$ is the fundamental solution of $\De^{2}$ in $\R^{n}$
and define
\begin{equation}
u_{2}=I_{2}(dP\wedge du),\label{eq: definition of u2}
\end{equation}
where $I_{2}$ is the fundamental solution of $-\De$ in $\R^{n}$.
It follows that
\begin{eqnarray*}
\De^{2}u_{11}+\De^{2}u_{12}=\De^{2}\tilde{u}_{1} & \text{and} & \De u_{2}=\De\tilde{u}_{2}
\end{eqnarray*}
on $B_{1}$. Set $h=d\tilde{u}_{1}-du_{11}-du_{12}+d^{\ast}\tilde{u}_{2}-d^{\ast}u_{2}+\tilde{h}$
so that
\[
\De^{2}h=0\qquad\text{in }B_{1}.
\]
We obtain the decomposition
\begin{eqnarray}
Pdu=du_{11}+du_{12}+d^{\ast}u_{2}+h &  & \text{in }B_{1}.\label{eq: modified hodge decomposition}
\end{eqnarray}

To obtain the Morrey decay estimates of $\na u$ and $\na^{2}u$, it suffices to estimate that of the
components $u_{11},u_{12}$ and $u_{2}$.

First we estimate $\na u_{11}$. From the definition \eqref{eq: definition of u11}
of $u_{11}$, it holds
\[
\na u_{11}=\na I_{4}\ast\left({\rm div}^{2}(D_{P}\otimes\na u)+{\rm div}(E_{P}\cdot\na u)+G_{P}\cdot\na u+\ast\D\De\xi\cdot Pdu+\De(\na P\na u)\right).
\]

Let $J_{1}=I_{4}\left({\rm div}^{2}(D_{P}\otimes\na u)+{\rm div}(E_{P}\cdot\na u)+\De(\na P\na u)\right)$. Then
\[
\na J_{1}\approx\na^{3}I_{4}(D_{P}\na u+\na P\na u)+\na^{2}I_{4}(E_{P}\na u),
\]
which implies that
\begin{equation}
|\na J_{1}|\lesssim I_{1}\left(|D_{P}||\na u|+|\na P||\na u|\right)+I_{2}\left(|E_{P}||\na u|\right).\label{eq: first order J1}
\end{equation}
Applying the growth condition \eqref{eq: growth conditions 2} gives
$$|D_{P}||\na u|+|\na P||\na u|\lesssim\left(|\na u|+|\na P|\right)|\na u|,$$
and
$$|E_{P}||\na u|\lesssim\left(\left|\nabla^{2}u\right|+|\nabla u|^{2}+\left|\nabla^{2}P\right|+|\nabla P|^{2}\right)|\na u|.$$
Since $\na P,\na u\in M_{\ast}^{4,n-4}(\R^{n})$ and $\na^{2}u,\na^{2}P\in M_{\ast}^{2,n-4}(\R^{n})$,
the H\"older inequality \eqref{eq: Morrey-Holder-2} implies that $D_{P}\na u\in M_{\ast}^{2,n-4}(\R^{n})$
and $E_{P}\na u\in M_{\ast}^{4/3,n-4}(\R^{n})$, together with estimates
\begin{equation}
\begin{aligned}\left\Vert D_{P}\na u\right\Vert _{M_{\ast}^{2,n-4}(\R^{n})} & \lesssim
\left(\|\na P\|_{M_{\ast}^{4,n-4}(\R^{n})}+\|\na u\|_{M_{\ast}^{4,n-4}(\R^{n})}\right)\|\na u\|_{M_{\ast}^{4,n-4}(\R^{n})}\\
&\lesssim\ep\|\na u\|_{M_{\ast}^{4,n-4}(\R^{n})},\end{aligned}
\label{eq: integrand of J-1-1}
\end{equation}
and
\begin{equation}
\begin{aligned}\left\Vert E_{P}\na u\right\Vert _{M_{\ast}^{4/3,n-4}(\R^{n})} & \lesssim\ep\|\na u\|_{M_{\ast}^{4,n-4}(\R^{n})}.
\end{aligned}
\label{eq: integrand of J-1-2}
\end{equation}
Here we used the bounded extension of $u,  P$ from  $M_{2, \ast}^{2,n-4}(B_1)$  into $M_{2,\ast}^{2,n-4}(\R^n)$ (see Proposition \ref{prop:Morrey-Sobolev extension}) and the smallness assumption
\eqref{eq: smallness assumption}. By Proposition \ref{prop: Ho},
\[
I_{1}\colon M_{\ast}^{2,n-4}(\R^{n})\to M_{\ast}^{4,n-4}(\R^{n})
\]
and
\[
I_{2}\colon M_{\ast}^{4/3,n-4}(\R^{n})\to M_{\ast}^{4,n-4}(\R^{n})
\]
are bounded operators. Thus from \eqref{eq: first order J1} and the
above estimates we deduce
\[
\|\na J_{1}\|_{M_{\ast}^{4,n-4}(\R^{n})}\lesssim\ep\|\na u\|_{M_{\ast}^{4,n-4}(\R^{n})}.
\]
Using the bounded extension $\|\na u\|_{M_{\ast}^{4,n-4}(\R^{n})}\lesssim\|\na u\|_{M_{\ast}^{4,n-4}(B_{1})}$,
it follows
\[
\|\na J_{1}\|_{M_{\ast}^{4,n-4}(B_{1})}\lesssim\ep\|\na u\|_{M_{\ast}^{4,n-4}(B_{1})}.
\]

Let $J_{2}=I_{4}(G_{P}\ast\na u)$. This is the most difficult term to estimate
and we need to exploit the full nonlinearity of $G_{P}$.
By \eqref{eq: Morrey-Holder-2} and the inequality \eqref{eq: Morrey-Holder},
and the fact $|\na u|,|\na P|\in M^{4,n-4}$, $|\na^{2}u|,|\na^{2}P|\in M^{2,n-4}$,
we infer that
\[
|G_{P}\na u|\lesssim\left(\left|\nabla^{2}u\right|+\left|\nabla^{2}P\right|\right)(|\nabla u|+|\nabla P|)|\na u|+\left(|\nabla u|^{3}+|\nabla P|^{3}\right)|\na u|\in M^{1,n-4}(\R^{n})
\]
 with estimates
\begin{equation}
\begin{aligned} & \left\Vert |G_{P}||\na u|\right\Vert _{M^{1,n-4}(\R^{n})}\\
 & \lesssim\|\na u\|_{M^{4,n-4}}\left(\|\na^{2}u\|_{M^{2,n-4}}+\|\na^{2}P\|_{M^{2,n-4}}\right)\left(\|\na u\|_{M^{4,n-4}}+\|\na P\|_{M^{4,n-4}}\right)\\
 & \quad+\|\na u\|_{M^{4,n-4}}\left(\|\na u\|_{M^{4,n-4}}^{3}+\|\na P\|_{M^{4,n-4}}^{3}\right).
\end{aligned}
\label{eq: Ingerand of J2}
\end{equation}
Combining the estimate \eqref{eq: estimate of gauge transform} of
$\na P$  with \eqref{eq: Ingerand of J2} yields
\[
\left\Vert |G_{P}||\na u|\right\Vert _{M^{1,n-4}(\R^{n})}\lesssim (\|\na^{2}u\|_{M^{2,n-4}}+\|\na u\|_{M^{4,n-4}}^{2})\|\na u\|_{M^{4,n-4}}^{2}\lesssim\|\na u\|_{M^{4,n-4}(\R^{n})}^{2}.
\]
Therefore, applying the bounded operator $I_{3}\colon M^{1,n-4}(\R^{n})\to M_{\ast}^{4,n-4}(\R^{n})$
by Proposition \ref{prop: Riesz in Morrey}, we arrive at
\[
\|\na J_{2}\|_{M_{\ast}^{4,n-4}(\R^{n})}\lesssim \left\Vert |G_{P}||\na u|\right\Vert _{M^{1,n-4}(\R^{n})}\lesssim \|\na u\|_{M^{4,n-4}(\R^{n})}^{2}.
\]
Thus we conclude
\[
\|\na J_{2}\|_{M_{\ast}^{4,n-4}(B_{1})}\lesssim\|\na u\|_{M^{4,n-4}(B_{1})}^{2}.
\]

Let $J_{3}=I_{4}(\ast\D\De\xi\cdot Pdu)$. Integrating by parts gives (up to signs)
$$J_{3}=\int\D\De\xi\wedge I_{4}Pdu=\int\De\xi\wedge(dI_{4}P+I_{4}dP)\wedge du.$$
Thus
\[
|\na J_{3}|\lesssim I_{2}\left(|\na^{2}\xi||\na u|\right)+I_{3}\left(|\na^{2}\xi||\na u||\na P|\right).
\]
As $|\na^{2}\xi||\na u|\in M_{\ast}^{2,n-4}\cdot M_{\ast}^{4,n-4}\subset M_{\ast}^{\frac{4}{3},n-4},$
we infer that $I_{2}\left(|\na^{2}\xi||\na u|\right)\in M_{\ast}^{4,n-4}$
as that of $J_{1}$ with estimate
\[
\left\Vert I_{2}\left(|\na^{2}\xi||\na u|\right)\right\Vert _{M_{\ast}^{4,n-4}(\R^{n})}\lesssim\ep\|\na u\|_{M_{\ast}^{4,n-4}(\R^{n})}.
\]
For the second term, we have $|\na^{2}\xi||\na u||\na P|\in M^{1,n-4}$.
As that of $J_{2}$, we obtain
\[
\left\Vert I_{3}\left(|\na^{2}\xi||\na u||\na P|\right)\right\Vert _{M_{\ast}^{4,n-4}(\R^{n})}\lesssim\|\na u\|_{M^{4,n-4}(\R^{n})}^{2}.
\]
Consequently,
\[
\left\Vert \na J_{3}\right\Vert _{M_{\ast}^{4,n-4}(\R^{n})}\lesssim\ep\|\na u\|_{M_{\ast}^{4,n-4}(\R^{n})}+\|\na u\|_{M^{4,n-4}(\R^{n})}^{2}.
\]
Using the bounded extension of $u$ gives
\[
\left\Vert \na J_{3}\right\Vert _{M_{\ast}^{4,n-4}(B_{1})}\lesssim\ep\|\na u\|_{M_{\ast}^{4,n-4}(B_{1})}+\|\na u\|_{M^{4,n-4}(B_{1})}^{2}.
\]

Taking the three estimates involving $\na J_{1},\na J_{2},\na J_{3}$,
we derive
\[
\left\Vert \na u_{11}\right\Vert _{M_{\ast}^{4,n-4}(B_{1})}\lesssim\ep\|\na u\|_{M_{\ast}^{4,n-4}(B_{1})}+\|\na u\|_{M^{4,n-4}(B_{1})}^{2}.
\]
Applying the inequality \eqref{Ineq: Struwe} and the smallness assumption \eqref{eq: smallness assumption}  and the embedding $M_{\ast}^{4,n-4}(B_{1})\subset M^{1,n-1}(B_{1})$, we find that
\[
\|\na u\|_{M^{4,n-4}(B_{1})}^{2}\lesssim\ep\|\na u\|_{M^{1,n-1}(B_{1})}\lesssim\ep\|\na u\|_{M_{\ast}^{4,n-4}(B_{1})}.
\]
Thus we obtain the estimate of $u_{11}$ as
\begin{equation}
\left\Vert \na u_{11}\right\Vert _{M_{\ast}^{4,n-4}(B_{1})}\lesssim\ep\|\na u\|_{M_{\ast}^{4,n-4}(B_{1})}.\label{eq: 1rd estimate of u11}
\end{equation}

The estimate of $u_{12}$ is standard.
Since $f\in M^{1,n-4+\alpha}$, $|\nabla u_{12}|\approx I_3(Pf)$ and standard
potential theory, Proposition \ref{prop: Riesz in Morrey}, gives $\nabla u_{12}\in M^{\frac{4-\al}{1-\al},n-4+\alpha}_\ast$. Notice that for $0<\al<1$, $\frac{4-\al}{1-\al}>4$ so we have
\begin{equation}
\|\na u_{12}\|_{M_{\ast}^{4,n-4}(B_{r})}\lesssim r^\alpha\|\nabla u_{12}\|_{M^{\frac{4-\al}{1-\al},n-4+\alpha}_\ast}\lesssim\|f\|_{M^{1,n-4+\alpha}(B_{1})}r^{\al}\label{eq: 1 order estimate of u12}
\end{equation}
 for any $r>0$. Here we have used the fact that $f\equiv0$ on $B_{1}^{c}$.

Combining the above estimates \eqref{eq: 1rd estimate of u11} and
\eqref{eq: 1 order estimate of u12}, we deduce that, for any $0<r\le1$,
\begin{equation}
\|\na u_{11}\|_{M_{\ast}^{4,n-4}(B_{r})}+\|\na u_{12}\|_{M_{\ast}^{4,n-4}(B_{r})}\lesssim\ep\|\na u\|_{M_{\ast}^{4,n-4}(B_{1})}+\|f\|_{M^{1,n-4+\alpha}(B_{1})}r^{\alpha}.\label{eq: estimate of u1}
\end{equation}

It remains to estimate $u_{2}$ and $h$. Since $u_{2}=I_{2}(dP\wedge du)$,
we have $|\na u_{2}|\lesssim I_{1}(|\na P||\na u|)$. As that of $J_{1}$,
we obtain
\begin{equation}
\|\na u_{2}\|_{M_{\ast}^{4,n-4}(B_{1})}\lesssim\ep\|\na u\|_{M_{\ast}^{4,n-4}(B_{1})}.\label{eq: estimate of u2}
\end{equation}
Since $h$ is biharmonic, for any $x\in B_{1}$ with $B_{2r}(x)\subset B_{1}$,
there holds
\[
\max_{B_{r}(x)}|\na h|\le C\fint_{B_{2r}(x)}|\na h|.
\]
So for any $x\in B_{1/2}$ and $0<r<1/2$,
\[
\|\na h\|_{L^{4,\wq}(B_{r}(x))}^{4}\le\int_{B_{r}(x)}|\na h|^{4}\lesssim r^{n}\max_{B_{1/2}}|\na h|^{4}\lesssim r^{n}\left(\fint_{B_{1}}|\na h|\right)^{4}\lesssim r^{n}\|\na h\|_{M_{\ast}^{4,n-4}(B_{1})}^{4}.
\]
That is,
\[
r^{\frac{n-4}{4}}\|\na h\|_{L^{4,\wq}(B_{r}(x))}\le r\|\na h\|_{M_{\ast}^{4,n-4}(B_{1})}.
\]
Hence
\begin{equation}
\|\na h\|_{M_{\ast}^{4,n-4}(B_{r})}=\sup_{x\in B_{r},0<s<2r}\left(s^{\frac{n-4}{4}}\|\na h\|_{L^{4,\wq}(B_{s}(x))}\right)\lesssim r\|\na h\|_{M_{\ast}^{4,n-4}(B_{1})}.\label{eq: estimate of h}
\end{equation}

Now we can obtain the decay estimate for $\na u$.
For any $0<\tau<1/2$, combining \eqref{eq: estimate of u1}, \eqref{eq: estimate of u2} and \eqref{eq: estimate of h}
gives
\[
\begin{aligned}\|\na &u\|_{M_{\ast}^{4,n-4}(B_{\tau})}\\
	 & \lesssim\|\na h\|_{M_{\ast}^{4,n-4}(B_{\tau})}+\|\na u_{11}\|_{M_{\ast}^{4,n-4}(B_{\tau})}+\|\na u_{12}\|_{M_{\ast}^{4,n-4}(B_{\tau})}+\|\na u_{2}\|_{M_{\ast}^{4,n-4}(B_{\tau})}\\
 & \lesssim\tau\|\na h\|_{M_{\ast}^{4,n-4}(B_{1})}+\ep\|\na u\|_{M_{\ast}^{4,n-4}(B_{1})}+\|f\|_{M^{1,n-4+\alpha}(B_{1})}\tau^{\al}\\
 & \lesssim\tau\Big(\|\na u\|_{M_{\ast}^{4,n-4}(B_{1})}+\|\na u_{11}\|_{M_{\ast}^{4,n-4}(B_{1})}+\|\na u_{12}\|_{M_{\ast}^{4,n-4}(B_{1})}\\
 & \qquad+\|\na u_{2}\|_{M_{\ast}^{4,n-4}(B_{1})}\Big)+\ep\|\na u\|_{M_{\ast}^{4,n-4}(B_{1})}+\tau^{\al}\|f\|_{M^{1,n-4+\alpha}(B_{1})}\\
 & \le C\left(\tau+\ep\right)\|\na u\|_{M_{\ast}^{4,n-4}(B_{1})}+C\tau^{\al}\|f\|_{M^{1,n-4+\alpha}(B_{1})}
\end{aligned}
\]
for some $C>0$ independent of $\tau$ and $\ep$. Recall that $0<\al<1$.
Take $\be\in(\al,1)$. Then take $\tau=r_{0}$ small enough such that
$2Cr_{0}<r_{0}^{\be}$, and then choose $\ep\le r_{0}$. We obtain
\[
\|\na u\|_{M_{\ast}^{4,n-4}(B_{r_{0}})}\le r_{0}^{\be}\|\na u\|_{M_{\ast}^{4,n-4}(B_{1})}+\|f\|_{M^{1,n-4+\alpha}(B_{1})}r_{0}^{\al}.
\]
Finally, using a standard scaling and translation and iteration argument,
there holds, for any $x\in B_{1/2}$ and $0<r<1$,
\[
\|\na u\|_{M_{\ast}^{4,n-4}(B_{r}(x))}\le Cr^{\al}\left(\|\na u\|_{M_{\ast}^{4,n-4}(B_{1})}+\|f\|_{M^{1,n-4+\alpha}(B_{1})}\right).
\]
In particular, this implies that for any $x\in B_{1/2}$ and $0<r<1$,
\[
\|\na u\|_{L^{4,\wq}(B_{r}(x))}^{4}\le Cr^{n-4+4\al}\left(\|\na u\|_{M_{\ast}^{4,n-4}(B_{1})}+\|f\|_{M^{1,n-4+\alpha}(B_{1})}\right)^{4}.
\]
Hence $\na u\in M_{\ast}^{4,n-4+4\al}(B_{1/2})$ and the desired estimate
\eqref{eq: 1 order decay estimate} follows.

Next we derive the decay of $\na^{2}u$ and the proof is similar to the one given above.
First estimate $\na^{2}u_{11}$. Using the same notations, we have
$$\na^{2}J_{1}\approx\na^{4}I_{4}(D_{P}\na u)+\na^{3}I_{4}(E_{P}\na u).$$
Since  $\na^{4}I_{4}$ is a singular integral operator,  Proposition \ref{prop: Adams-singular integral operator} implies that
\[
\na^{4}I_{4}\colon M^{p,\la}(\R^{n})\to M^{p,\la}(\R^{n})
\]
is a bounded operator. Thus
\[
\|\na^{4}I_{4}(D_{P}\na u)\|_{M^{2,n-4}}\lesssim\left(\|\na P\|_{M^{4,n-4}}+\|\na u\|_{M^{4,n-4}}\right)\|\na u\|_{M^{4,n-4}}\lesssim\|\na u\|_{M^{4,n-4}(\R^{n})}^{2},
\]
where the second inequality follows from inequality \eqref{eq: estimate of gauge transform}.
Using the embedding $M^{2,n-4}\subset M_{\ast}^{2,n-4}$, the inequality
\eqref{Ineq: Struwe} and the smallness assumption \eqref{eq: smallness assumption}
of $u$ as before, we deduce
\[
\|\na^{4}I_{4}(D_{P}\na u)\|_{M_{\ast}^{2,n-4}(B_{1})}\lesssim\ep\|\na u\|_{M_{\ast}^{4,n-4}(B_{1})}.
\]
For the second term, combining \eqref{eq: integrand of J-1-2} and
the boundedness of $$I_{1}\colon M_{\ast}^{4/3,n-4}(\R^{n})\to M_{\ast}^{2,n-4}(\R^{n})$$
by Proposition \ref{prop: Ho}, we infer
\[
\left\Vert \na^{3}I_{4}(E_{P}\na u)\right\Vert _{M_{\ast}^{2,n-4}(B_{1})}\lesssim\ep\|\na u\|_{M_{\ast}^{4,n-4}(B_{1})}.
\]
Hence
\[
\left\Vert \na^{2}J_{1}\right\Vert _{M_{\ast}^{2,n-4}(B_{1})}\lesssim\ep\|\na u\|_{M_{\ast}^{4,n-4}(B_{1})}.
\]

For $J_{2}$, we have
\[
|\na^{2}J_{2}|\lesssim I_{2}(|G_{P}||\na u|).
\]
Recall that $G_{P}\na u\in M^{1,n-4}(\R^{n})$ and estimate \eqref{eq: Ingerand of J2}
holds. Hence $\na^{2}J_{2}\in M_{\ast}^{2,n-4}(\R^{n})$ by Proposition
\ref{prop: Riesz in Morrey} with estimate
\[
\left\Vert \na^{2}J_{2}\right\Vert _{M_{\ast}^{2,n-4}}\lesssim\left\Vert |G_{P}||\na u|\right\Vert _{M^{1,n-4}}\lesssim\|\na u\|_{M^{4,n-4}(\R^{n})}^{2}.
\]
Again, applying inequality \eqref{Ineq: Struwe} yields
\[
\left\Vert \na^{2}J_{2}\right\Vert _{M_{\ast}^{2,n-4}(B_{1})}\lesssim\ep\|\na u\|_{M_{\ast}^{4,n-4}(B_{1})}.
\]

For $J_{3}$, we have
\[
|\na^{2}J_{3}|\lesssim I_{1}\left(|\na^{2}\xi||\na u|\right)+I_{2}\left(|\na^{2}\xi||\na u||\na P|\right).
\]
Similar to $J_{1}$ and $J_{2}$, we derive
\[
\left\Vert \na^{2}J_{3}\right\Vert _{M_{\ast}^{2,n-4}(B_{1})}\lesssim\ep\|\na u\|_{M_{\ast}^{4,n-4}(B_{1})}.
\]

All together we conclude that
\begin{equation}
\left\Vert \na^{2}u_{11}\right\Vert _{M_{\ast}^{2,n-4}(B_{1})}\lesssim\ep\|\na u\|_{M_{\ast}^{4,n-4}(B_{1})}.\label{eq: 2nd order estimate of u11}
\end{equation}

For $u_{12}$, since $f\in M^{1,n-4+\alpha}$ and $|\nabla^2 u_{12}|\approx I_2(Pf)$,   Proposition \ref{prop: Riesz in Morrey} gives $\nabla^2 u_{12}\in M^{\frac{4-\al}{2-\al},n-4+\alpha}_\ast$. Notice that for $0<\al<1$, $\frac{4-\al}{2-\al}>2$. So similar to \eqref{eq: 1 order estimate of u12}, we obtain for any
$0<r<\wq$,
\begin{equation}
\|\na^{2}u_{12}\|_{M_{\ast}^{2,n-4}(B_{r})}\lesssim\|f\|_{M^{1,n-4+\alpha}(B_{1})}r^{\al}.\label{eq: 2nd order estimate of u12}
\end{equation}

For the term $u_{2}$, we have $|\na^{2}u_{2}|\lesssim I_{0}(|\na P||\na u|)\in M^{2,n-4}$
with
\begin{equation}
\|\na^{2}u_{2}\|_{M_{\ast}^{2,n-4}}\le\|\na^{2}g\|_{M^{2,n-4}}\lesssim\|\na u\|_{M^{4,n-4}}^{2}\lesssim\ep\|\na u\|_{M_{\ast}^{4,n-4}(B_{1})}.\label{eq: 2nd order estimate of u2}
\end{equation}
Similarly dispose the biharmonic 1-form $h$, for any $0<r<1$, there
holds
\begin{equation}
\|\na^{2}h\|_{M_{\ast}^{2,n-4}(B_{r})}\lesssim r\|\na^{2}h\|_{M_{\ast}^{2,n-4}(B_{1})}.\label{eq: 2nd order estimate of h}
\end{equation}

Combining estimates \eqref{eq: 2nd order estimate of u11},
\eqref{eq: 2nd order estimate of u12}, \eqref{eq: 2nd order estimate of u2}
and \eqref{eq: 2nd order estimate of h} yields, for any $0<r<1$,
\[
\begin{aligned}
& \|\na^{2}u\|_{M_{\ast}^{2,n-4}(B_{r})}+\|\na u\|_{M_{\ast}^{4,n-4}(B_{r})}\\
\lesssim & (r+\ep)\left(\|\na^{2}u\|_{M_{\ast}^{2,n-4}(B_{1})}+\|\na u\|_{M_{\ast}^{4,n-4}(B_{1})}\right)+\|f\|_{M^{1,n-4+\alpha}(B_{1})}r^{\al}.
\end{aligned}\]
Similar iteration, scaling and translation arguments give
\[
\begin{aligned}
&\|\na^{2}u\|_{M_{\ast}^{2,n-4}(B_{r})}+\|\na u\|_{M_{\ast}^{4,n-4}(B_{r})}\\
\lesssim &\left(\|\na^{2}u\|_{M_{\ast}^{2,n-4}(B_{1})}+\|\na u\|_{M_{\ast}^{4,n-4}(B_{1})}+\|f\|_{M^{1,n-4+\alpha}(B_{1})}\right)r^{\al}.
\end{aligned}
\]
The proof is complete.
\end{proof}

\begin{proof}[Proof of Corollary \ref{coro:Holder}]
	It follows from Theorem \ref{thm: Morrey decay estimate} and Proposition \ref{prop: Morrey's DGT}.
\end{proof}

\section{$L^{p}$ regularity theory}

In this section we prove Theorem \ref{thm: main result} and Theorem  \ref{thm:p small}. We will write
\[p_1=\frac {np}{n-p},\quad p_2=\frac {np}{n-2p},\quad p_3=\frac {np}{n-3p},\]
whenever these are positive numbers. For $p<n$, set
\begin{equation}\label{eq:alpha}
	\al=4-n/p.
\end{equation}


Roughly speaking, Theorem \ref{thm: main result} and Theorem  \ref{thm:p small} follow from the Morrey estimate of the previous section and an iteration argument. Along the iteration the constant $\ep$ should become smaller and smaller. Fortunately, the iteration stops after finitely many steps. Thus we  can always choose a sufficiently small  $\ep$  in the very beginning such that  the whole iteration proceeds. As in the previous proofs, the Gauge transform plays a central role.

\subsection{Case 1: $n/4<p<n/3$}

In this subsection we prove Theorem  \ref{thm: main result}
in the case $n/4<p<n/3$. Recall that our initial assumption is that
$\na u\in M^{4,n-4}(B_{1}),\na^{2}u\in M^{2,n-4}(B_{1})$ hold with
the smallness assumption \eqref{eq: smallness assumption}.
Thus we can choose $\ep$ sufficiently small such that we have the improvement
\[
\na u\in M_{\ast}^{4,n-4+4\al}(B_{1/2})\quad\text{and}\quad\na^{2}u\in M_{\ast}^{2,n-4+2\al}(B_{1/2}),
\]
where $\al=4-n/p\in(0,1)$.
At this moment, due to the strong nonlinearity, the regularity
of the function
\[
|G_{P}\na u|\lesssim\left(\left|\nabla^{2}u\right|+\left|\nabla^{2}P\right|\right)(|\nabla u|+|\nabla P|)|\na u|+\left(|\nabla u|^{3}+|\nabla P|^{3}\right)|\na u|
\]
will be too weak to iterate.

Fortunately we have the following two observations. The first one is that the second order weak Morrey regularity implies:
\begin{equation}
\na u\in L^{2\chi}\cap M_{\ast}^{2\chi,n-2\chi(1-\al)}(B_{1/4})\label{eq: 1st order improvement of Lebesgue}
\end{equation}
where
\begin{equation}
\chi\equiv(2-\al)/(1-\al)>2.\label{eq: def of chi}
\end{equation}

To find this, select $\eta\in C_{0}^{\wq}(B_{1/2})$ with $\eta\equiv1$
on $B_{1/4}$. An elementary calculation shows that $\na(u\eta)\in M_{\ast}^{4,n-4(1-\al)}(B_{1/2})$
and $\na^{2}(\eta u)\in M_{\ast}^{2,n-2(2-\al)}(B_{1/2})$. Set $\eta u\equiv0$
outside $B_{1/2}$.  \eqref{eq: Ho+Adams} implies that
\[
\na(\eta u)=\na I_{2}(-\De(\eta u))\approx I_{1}(\De u)\in L^{2\chi}\cap M_{\ast}^{2\chi,n-2\chi(1-\al)}(\R^{n})
\]
with estimates
\[
\|\na(\eta u)\|_{L^{2\chi}}+\|\na(\eta u)\|_{M_{\ast}^{2\chi,n-2\chi(1-\al)}}\lesssim \|\De(\eta u)\|_{L^{2}}+\|\De(\eta u)\|_{M_{\ast}^{2,n-2(2-\al)}}.
\]
This yields \eqref{eq: 1st order improvement of Lebesgue} for $\na u$
with
\[
\|\na u\|_{L^{2\chi}(B_{1/4})}+\|\na u\|_{M_{\ast}^{2\chi,n-2\chi(1-\al)}(B_{1/4})}\lesssim \|\na u\|_{M_{\ast}^{4,n-4(1-\al)}(B_{1/2})}+\|\na^{2}u\|_{M_{\ast}^{2,n-2(2-\al)}(B_{1/2})}.
\]

The second observation is:

\begin{lemma}\label{lem: Gauge transform-general-version}  There exist $\ep=\ep(n,m)>0$
and $C=C(n,m)>0$ with the following property: For every $\Om\in M_{1}^{2,n-4}\cap M^{4,n-4}(B_{1/2},so_{m}\otimes\wedge^{1}\R^{n})$
with
\[
\|\na\Om\|_{M^{2,n-4}(B_{1/2})}+\|\Om\|_{M^{4,n-4}(B_{1/2})}\le\ep,
\]
there exist $P\in W^{2,2}(B_{1/2},SO_{m})$ and $\xi\in W^{2,2}(B_{1/2},so_{m}\otimes\wedge^{n-2}\R^{n})$
such that Lemma \ref{lem: Gauge transform} holds on $B_{1/2}$.

In addition, if $\Om\in M_{\ast}^{4,n-4+4\al}(B_{1/2})$ and $\na\Om\in M_{\ast}^{2,n-4+2\al}(B_{1/2})$,
then we further have $\na P,\na\xi\in M_{\ast}^{4,n-4+4\al}(B_{1/2})$,
$\na^{2}P,\na^{2}\xi\in M_{\ast}^{2,n-4+2\al}(B_{1/2})$ together
with
\begin{equation}
\|\na P\|_{M_{\ast}^{4,n-4+4\al}(B_{1/2})}+\|\na\xi\|_{M_{\ast}^{4,n-4+4\al}(B_{1/2})}\le C\|\Om\|_{M_{\ast}^{4,n-4+4\al}(B_{1/2})},\label{eq: improved gauge regularity 1}
\end{equation}
and
\begin{equation}
\begin{aligned}
  &\|\na^{2}P\|_{M_{\ast}^{2,n-4+2\al}(B_{1/2})}+\|\na^{2}\xi\|_{M_{\ast}^{2,n-4+2\al}(B_{1/2})}\\
  \le & C\left(\|\na\Om\|_{M_{\ast}^{2,n-4+2\al}(B_{1/2})}
+\|\Om\|_{M_{\ast}^{4,n-4+4\al}(B_{1/2})}\right).
\end{aligned}
\label{eq: improved gauge regularity 2}
\end{equation}
 \end{lemma}
\begin{proof}
The existence of $P,\xi$ follows from the same method as that of Lemma \ref{lem: Gauge transform}. For the proof of estimates \eqref{eq: improved gauge regularity 1} and \eqref{eq: improved gauge regularity 2}, see Lemma \ref{lem: higher order apriori estimate for gauge} in the Appendix.
\end{proof}
Let $P,\xi$ be obtained as in Lemma \ref{lem: Gauge transform-general-version}. By the first observation, we
have
\begin{equation}
\na P,\na\xi\in L^{2\chi}\cap M_{\ast}^{2\chi,n-2\chi(1-\al)}(B_{1/4})\label{eq: 1st order improvement of p-xi}
\end{equation}
and
\[
\|\na P,\na\xi\|_{L^{2\chi}\cap M_{\ast}^{2\chi,n-2\chi(1-\al)}(B_{1/4})}\lesssim\|\na\Om\|_{M_{\ast}^{2,n-4+2\al}(B_{1/2})}+\|\Om\|_{M_{\ast}^{4,n-4+4\al}(B_{1/2})}.
\]
Thus we deduce from the growth assumption on $\Om$ that
\[
\|\na P,\na\xi\|_{L^{2\chi}\cap M_{\ast}^{2\chi,n-2\chi(1-\al)}(B_{1/4})}\lesssim\|\na^{2}u\|_{M_{\ast}^{2,n-4+2\al}(B_{1/2})}+\|\na u\|_{M_{\ast}^{4,n-4+4\al}(B_{1/2})}.
\]

We transform the system \eqref{eq: Struwe system} on $B_{1/4}$
to obtain the gauge equivalent system \eqref{eq: (35)}. Then we
extend all the functions from $B_{1/4}$ into $\R^{n}$ with compact
supports in $B_{2}$ in a bounded way, and define similarly  $u_{11},u_{12}$,
$u_{2}$ and a biharmonic 1-form $h$ on $B_{1/4}$ as that of \eqref{eq: definition of u11}, \eqref{eq: def of u12} and \eqref{eq: definition of u2} such that
$Pdu=du_{11}+du_{12}+d^{\ast}u_{2}+h$ on $B_{1/4}$.

Our aim is to improve the regularity of $\na^{2}u$ through the gauge
equivalent system \eqref{eq: (35)}.

\begin{claim}\label{claim: 1st iteration} Let $p_{2}=\frac{np}{n-2p}$
and $\chi$ be defined as in \eqref{eq: def of chi}. Then
\begin{equation}
\na^{2}u\in\begin{cases}
L^{\chi}\cap M_{\ast}^{\chi,n-\chi(2-\al)}(B_{\frac{1}{4}}) & \text{if }\chi<p_{2},\\
L^{p_{2}}(B_{\frac{1}{4}}) & \text{if }\chi\ge p_{2},
\end{cases}\label{eq: 1 time iteration}
\end{equation}
\end{claim}
\begin{proof}
Hereafter all the norms are taken on the whole space $\R^{n}$ unless specified.  We first deduce the regularity of $\na^{2}u_{11}$.

For the first term $J_{1}$, \eqref{eq: 1st order improvement of Lebesgue},
\eqref{eq: 1st order improvement of p-xi} and H\"older's inequality
\eqref{eq: Morrey-Holder-2} imply
\[
\left(|\na^{2}u|+|\na^{2}P|+|\na u|^{2}+|\na P|^{2}\right)(|\nabla u|+|\nabla P|)\in M_{\ast}^{\frac{2\chi}{\chi+1},n-4+2\al}\cap L^{\frac{2\chi}{\chi+1}}.
\]
Since
\[
\na^{2}J_{1}\approx I_{1}\left(\left(|\na^{2}u|+|\na^{2}P|+|\na u|^{2}+|\na P|^{2}\right)(|\nabla u|+|\nabla P|)\right),
\]
and by \eqref{eq: Ho+Adams}
\[
I_{1}\colon L^{\frac{2\chi}{\chi+1}}\cap M_{\ast}^{\frac{2\chi}{\chi+1},n-4+2\al}\to L^{\chi}\cap M_{\ast}^{\chi,n-4+2\al}
\]
is a bounded operator, we obtain $\na^{2}J_{1}\in M_{\ast}^{\chi,n-4+2\al}(\R^{n})\cap L^{\chi}(\R^{n})$ with
\[
\left\Vert \na^{2}J_{1}\right\Vert _{M_{\ast}^{\chi,n-4+2\al}\cap L^{\chi}(\R^{n})}\lesssim\|\na u\|_{L^{\chi}\cap M_{\ast}^{\chi,n-\chi(2-\al)}(B_{1/4})}\|\na^{2}u\|_{L^{2}\cap M_{\ast}^{2,n-2(2-\al)}(B_{1/4})}.
\]
By the weak Morrey estimate,
\[
\|\na u\|_{L^{\chi}\cap M_{\ast}^{\chi,n-\chi(2-\al)}(B_{1/4})}+\|\na^{2}u\|_{L^{2}\cap M_{\ast}^{2,n-2(2-\al)}(B_{1/4})}\lesssim\ep+\|f\|_{L^p(B_1)}.
\]
This in turn leads
\begin{equation}
\left\Vert \na^{2}J_{1}\right\Vert _{M_{\ast}^{\chi,n-4+2\al}\cap L^{\chi}(\R^{n})}\lesssim(\ep+\|f\|_{L^p(B_1)})\|\na^{2}u\|_{L^{2}\cap M_{\ast}^{2,n-2(2-\al)}(B_{1/4})}.\label{eq: nonlinear apriori esti-1}
\end{equation}

For the second term, we have $|\na^{2}J_{2}|\lesssim I_{2}(|G_{P}||\na u|)$
and
\[
|G_{P}|\le\left(\left|\nabla^{2}u\right|+\left|\nabla^{2}P\right|\right)(|\nabla u|+|\nabla P|)+\left(|\nabla u|^{3}+|\nabla P|^{3}\right).
\]
Recall that $\na u,\na P\in M_{\ast}^{2\chi,n-4+2\al}\cap L^{2\chi}$
and $\na^{2}u,\na^{2}P\in M_{\ast}^{2,n-4+2\al}\cap L^{2}$. So
\[
\left(|\na^{2}u|+|\na^{2}P|\right)(|\nabla u|+|\nabla P|)^{2}\in M_{\ast}^{\frac{2\chi}{\chi+2},n-4+2\al}\cap L^{\frac{2\chi}{\chi+2}},
\]
\[
(|\nabla u|+|\nabla P|)^{4}=(|\nabla u|+|\nabla P|)^{2}(|\nabla u|+|\nabla P|)^{2}\in M_{\ast}^{\frac{2\chi}{\chi+2},n-4+2\al}\cap L^{\frac{2\chi}{\chi+2}}.
\]
Here the first term can be regarded in the space $M_{\ast}^{4,n-4+2\al}$
in view of the embedding $M_{\ast}^{2\chi,n-4+2\al}(B_{1/2})\subset M_{\ast}^{4,n-4+2\al}(B_{1/2}).$
Thus
\[
G_{P}\na u\in M_{\ast}^{\frac{2\chi}{\chi+2},n-4+2\al}\cap L^{\frac{2\chi}{\chi+2}}.
\]
In the case $\al<2/3$, we may apply \eqref{eq: Ho+Adams} to deduce
the boundedness of
\[
I_{2}\colon M_{\ast}^{\frac{2\chi}{\chi+2},n-4+2\al}\cap L^{\frac{2\chi}{\chi+2}}\to M_{\ast}^{\frac{2(2-\al)}{2-3\al},n-4+2\al}\cap L^{\frac{2(2-\al)}{2-3\al}},
\]
which implies $\na^{2}J_{2}\in M_{\ast}^{\frac{2(2-\al)}{2-3\al},n-4+2\al}\cap L^{\frac{2(2-\al)}{2-3\al}}$.
Similar to
(\ref{eq: nonlinear apriori esti-1}), we can obtain
\begin{equation}
\left\Vert \na^{2}J_{2}\right\Vert _{M_{\ast}^{\frac{2(2-\al)}{2-3\al},n-4+2\al}\cap L^{\frac{2(2-\al)}{2-3\al}}}\lesssim(\ep+\|f\|_{L^p(B_{1}}))^{a}\|\na u\|_{L^{4}\cap M_{\ast}^{4,n-4(1-\al)}(B_{1/4})}.\label{eq: nonlinear apriori esti-2}
\end{equation}
 for some $a>0$.

For the third term $|\na^{2}J_{3}|\lesssim I_{1}(|\na^{2}\xi||\na u|)+I_{2}(|\na^{2}\xi||\na u|||\na P|)$,
the same estimates as that of $J_{1}$ and $J_{2}$ imply
\[
I_{1}(|\na^{2}\xi||\na u|)\in M_{\ast}^{\chi,n-4+2\al}(\R^{n})\cap L^{\chi}(\R^{n})
\]
and when $\al<2/3$
\[
I_{2}(|\na^{2}\xi||\na u|||\na P|)\in M_{\ast}^{\frac{2(2-\al)}{2-3\al},n-4+2\al}\cap L^{\frac{2(2-\al)}{2-3\al}}.
\]
Note that if $\al<2/3$, then
\[
\frac{2(2-\al)}{2-3\al}=\frac{2-\al}{1-\frac{3}{2}\al}>\chi,
\]
and if $\al>\frac{2}{3}$, then the regularity of $\nabla^2 J_{i}$, $i=2,3$, become even better. All together, we may  conclude
\[
\na^{2}u_{11}\in M_{\ast}^{\chi,n-4+2\al}(B_{1/4})\cap L^{\chi}(B_{1/4}).
\]

Since $u_{12}\in W^{4,p}(\R^{n})$, $\na^{2}u_{12}\in W^{2,p}(\R^{n})\subset L^{p_{2}}$.
In particular, for any $s>0$,
\[
\int_{B_{s}(x)}|\na^{2}u_{12}|^{\chi}\lesssim\|\na^{2}J_{5}\|_{L^{p_{2}}}^{\chi}s^{n-2\chi+\chi\al}\lesssim\|f\|_{L^{p}}^{\chi}s^{n-2\chi+\chi\al}.
\]
That is
\[
\na^{2}u_{12}\in L^{p_{2}}\cap M^{\chi,n-2\chi+\chi\al}(\R^{n}).
\]

Similar to the estimate of $J_{1}$, one deduces
\[
\na^{2}u_{2}\in M_{\ast}^{\chi,n-4+2\al}(B_{1/4})\cap L^{\chi}(B_{1/4}).
\]
Note that the biharmonic 1-form $h$ is always smooth. Hence Claim
\ref{claim: 1st iteration} holds if $\chi\ge p_{2}$. In the case
$\chi<p_{2}=np/(n-2p)=n/(2-\al)$, observe that $n-4+2\al=n-2\chi+2\chi\al$.
So, for any $w\in M_{\ast}^{\chi,n-2\chi+2\chi\al}(B_{1/4})$ and
any $0<r<1/2$,
\[
\|w\|_{L^{\chi,\wq}(B_{r}(x))}^{\chi}\le\|w\|_{M_{\ast}^{\chi,n-2\chi+2\chi\al}(B_{1/4})}^{\chi}r^{n-2\chi+2\chi\al}\le\|w\|_{M_{\ast}^{\chi,n-2\chi+2\chi\al}(B_{1/4})}^{\chi}r^{n-2\chi+\chi\al}.
\]
That is,
\[
M_{\ast}^{\chi,n-4+2\al}(B_{1/4})\subset M^{\chi,n-2\chi+\chi\al}(B_{1/4}).
\]
Therefore,
\[
\na^{2}u_{11,}\na^{2}u_{2}\in L^{\chi}\cap M_{\ast}^{\chi,n-2\chi+\chi\al}(B_{1/4}).
\]
The proof of Claim \ref{claim: 1st iteration} is complete.
\end{proof}

Next we use  iteration to derive the optimal regularity of $\na u$ and $\na^{2}u$. 

\begin{claim}\label{claim: iteration} (Iteration lemma) Let $\chi=\la_{1}\le\la_{n}<p_{2}$
and set
\[
\la_{n+1}=\frac{\chi}{2}\la_{n}.
\]
If
\[
\na u\in L^{2\la_{n}}\cap M_{\ast}^{2\la_{n},n-2\la_{n}(1-\al)}(B_{1})\quad\text{and}\quad \na^{2}u\in L^{\la_{n}}\cap M_{\ast}^{\la_{n},n-\la_{n}(2-\al)}(B_{1})
\]
with sufficiently small norms and if $$\la_{n+1}<\chi,$$ then
\[
\na u\in L^{2\la_{n+1}}\cap M_{\ast}^{2\la_{n+1},n-2\la_{n+1}(1-\al)}(B_{1/2})\quad\text{and}\quad\na^{2}u\in L^{\la_{n+1}}\cap M_{\ast}^{\la_{n+1},n-\la_{n+1}(2-\al)}(B_{1/2}).
\]
\end{claim}
\begin{proof}
The improvement of $\na u$ on $B_{1/2}$ follows as before. We also
have the same regularity of $P,\xi$ as that of $u$ by the same arguments
as above. So we only need to deduce the regularity of $\na^{2}u$ and the arguments will be similar as in the previous step.

For the first term, we have
\[
\na^{2}J_{1}\approx I_{1}\left(\left(|\na^{2}u|+|\na^{2}P|\right)|\na u|+(|\nabla u|+|\nabla P|)|\na^{2}u|\right).
\]
H\"older's inequality gives
\[
\left(|\na^{2}u|+|\na^{2}P|\right)|\na u|+(|\nabla u|+|\nabla P|)|\na^{2}u|\in L^{\frac{\chi}{\chi+1}\la_{n}}\bigcap M_{\ast}^{\frac{\chi}{\chi+1}\la_{n},n-\frac{\chi}{\chi+1}\la_{n}(3-2\al)}
\]
and \eqref{eq: Ho+Adams} gives the boundedness of
\[
I_{1}\colon L^{\frac{\chi}{\chi+1}\la_{n}}\bigcap M_{\ast}^{\frac{\chi}{\chi+1}\la_{n},n-\frac{\chi}{\chi+1}\la_{n}(3-2\al)}\to L^{\frac{3-2\al}{2-2\al}\frac{\chi}{\chi+1}\la_{n}}\cap M_{\ast}^{\frac{3-2\al}{2-2\al}\frac{\chi}{\chi+1}\la_{n},n-\frac{\chi}{\chi+1}\la_{n}(3-2\al)}.
\]
Note that
\[
\frac{3-2\al}{2-2\al}\frac{\chi}{\chi+1}\la_{n}=\frac{\chi}{2}\la_{n}=\la_{n+1}.
\]
Hence
\[
\na^{2}J_{1}\in L^{\la_{n+1}}\cap M_{\ast}^{\la_{n+1},n-\la_{n+1}(2-2\al)}.
\]

For the second term, we have $|\na^{2}J_{2}|\lesssim I_{2}(|G_{P}||\na u|)$
and
\[
|G_{P}\na u|\le\left(\left|\nabla^{2}u\right|+\left|\nabla^{2}P\right|\right)(|\nabla u|+|\nabla P|)^{2}+\left(|\nabla u|+|\nabla P|\right)^{4}.
\]
Note that $2\la_{n+1}=\chi\la_{n}$ and so
\[
\left(|\na^{2}u|+|\na^{2}P|\right)(|\nabla u|+|\nabla P|)^{2}\in L^{\frac{\chi}{\chi+2}\la_{n}}\bigcap M_{\ast}^{\frac{\chi}{\chi+2}\la_{n},n-\frac{\chi}{\chi+2}\la_{n}(4-3\al)},
\]
\[
(|\nabla u|+|\nabla P|)^{4}\in L^{\frac{\chi}{4}\la_{n}}\bigcap M_{\ast}^{\frac{\chi}{4}\la_{n},n-\chi\la_{n}(1-\al)}.
\]
The second term has better regularity than the first one. In the case $\al<2/3$,
applying \eqref{eq: Ho+Adams} gives boundedness of
\[
I_{2}\colon L^{\frac{\chi}{\chi+2}\la_{n}}\bigcap M_{\ast}^{\frac{\chi}{\chi+2}\la_{n},n-\frac{\chi}{\chi+2}\la_{n}(4-3\al)}\to L^{\frac{4-3\al}{2-3\al}\frac{\chi}{\chi+2}\la_{n}}\bigcap M_{\ast}^{\frac{4-3\al}{2-3\al}\frac{\chi}{\chi+2}\la_{n},n-\frac{\chi}{\chi+2}\la_{n}(4-3\al)}.
\]
So
\[
\na^{2}J_{2}\in L^{\frac{4-3\al}{2-3\al}\frac{\chi}{\chi+2}\la_{n}}\bigcap M_{\ast}^{\frac{4-3\al}{2-3\al}\frac{\chi}{\chi+2}\la_{n},n-\frac{\chi}{\chi+2}\la_{n}(4-3\al)}.
\]

The third term can be splited into a sum of two terms with the same regularity
as that of $J_{1}$ and $J_{2}$. Since $\chi+2=(4-3\al)/(1-\al)$, we
have
\[
\tilde{\la}_{n+1}\equiv\frac{4-3\al}{2-3\al}\frac{\chi}{\chi+2}\la_{n}=\frac{2-\al}{2-3\al}\la_{n}=\frac{2(1-\al)}{2-3\al}\la_{n+1}.
\]
Hence $\tilde{\la}_{n+1}>\la_{n+1}$ and
\[
n-\tilde{\la}_{n+1}(2-3\al)=n-\la_{n+1}(2-2\al)>n-\la_{n+1}(2-\al).
\]
This implies
\[
L^{\tilde{\la}_{n+1}}\cap M_{\ast}^{\tilde{\la}_{n+1},n-\tilde{\la}_{n+1}(2-3\al)}(B_{1/2})\subset L^{\la_{n+1}}\cap M_{\ast}^{\la_{n+1},n-\la_{n+1}(2-\al)}(B_{1/2}).
\]
Consequently, we obtain
\[
\na^{2}u_{11}\in L^{\la_{n+1}}\cap M_{\ast}^{\la_{n+1},n-\la_{n+1}(2-\al)}(B_{1/2}).
\]

Note also that $\na^{2}u_{12}\in L^{p_{2}}.$ Thus if $\la_{n+1}<p_{2}$,
then
\[
\na^{2}u_{12}\in L^{\la_{n+1}}\cap M^{\la_{n+1},n-\la_{n+1}(2-\al)}(B_{1/2}).
\]

Similarly, we can deduce the result for $u_{2}$ and the biharmonic
part $h$. The proof of Claim \ref{claim: iteration} is complete.
\end{proof}

Since $\chi>2$, Claims \ref{claim: 1st iteration} and \ref{claim: iteration}
imply that after finitely many steps, this iteration will stop, whence $\na^{2}u\in L_{\loc}^{p_{2}}(B_{1})$. This in return implies
that $\na u\in L_{\loc}^{p_{3}}(B_{1})$ by the Sobolev embedding
theorem.

Now we can deduce the third order regularity of $u$. Rewrite the system \eqref{eq: Struwe system} as
\begin{equation}
\De^{2}u={\rm div}\left(I\right)+II,\label{eq: original equation rewritten}
\end{equation}
where
\[
I=D\na^{2}u+\na D\cdot\na u+E\cdot\na u+\na\Om\cdot\na u,
\]
\[
II=-\na\Om\cdot\De u+G\cdot\na u+f.
\]
By the growth assumption \eqref{eq:GC fourth order},
we know
\begin{equation}\label{eq: growth of I}
|I|\le C(|\na^{2}u|+|\na u|^{2})|\na u|,
\end{equation}
\[
|II|\le C\left(|\na^{2}u|^{2}+|\na u|^{2}|\na^{2}u|+|\na u|^{4}\right)+f.
\]
Since we have proved that $u\in W^{2,\frac{np}{n-2p}}$ and $n\le 4p$,
it follows that $2p\le p_{2}\le p_{3}/2$. Hence $I\in L_{\loc}^{\frac{np}{2n-5p}}\subset L_{\loc}^{\frac{np}{n-p}}$
and $II\in L_{\loc}^{p}$. Here the least regular term of $I$ and
$II$ are $\na^{2}u\na u$ and $f$, respectively.

Set $\De^{2}u_{1}={\rm div}(I)$
and $\De^{2}u_{2}=II$ in $B_{1}$. Standard elliptic regularity theory
implies $u_{1}\in W^{3,\frac{np}{n-p}}_{\loc}$ and $u_{2}\in W_{\loc}^{4,p}\subset W^{3,\frac{np}{n-p}}$.
As $u-u_{1}-u_{2}$ is a biharmonic function, we infer that
\[
u\in W_{\loc}^{3,\frac{np}{n-p}}(B_{1}).
\]

Next we derive the apriori estimate of $u$. By the Hodge decomposition,
we have the biharmonic 1-form $h$ satisfying
\[
h=Pdu-du_{11}-du_{12}-d^{\ast}u_{2}\qquad\text{in }B_{1}.
\]
By the Morrey estimates (Theorem \ref{thm: Morrey decay estimate}), we have
\[
\|h\|_{M_{\ast}^{4,n-4+4\al}(B_{3/4})}\lesssim\ep+\|f\|_{L^{p}(B_{1})}.
\]
In particular, this implies that
$\|h\|_{L^{1}(B_{3/4})}\lesssim\ep+\|f\|_{L^{p}(B_{1})}.$
Since $h$ is biharmonic, we infer that
\[
\|h\|_{L^{p_{3}}(B_{1/2})}\lesssim\|h\|_{L^{1}(B_{3/4})}\lesssim\ep+\|f\|_{L^{p}(B_{1})}.
\]
Returning to the Hodge decomposition, we have
\[
\|\na u\|_{L^{p_{3}}(B_{1/2})}\lesssim\|h\|_{L^{p_{3}}(B_{1/2})}
+\|\nabla u_{11}\|_{L^{p_{3}}(B_{1/2})}
+\|\nabla u_{12}\|_{L^{p_{3}}(B_{1/2})}+\|\nabla u_{2}\|_{L^{p_{3}}(B_{1/2})}.
\]
Using the potential theory, we can similarly estimate $\|\nabla u_{11}\|_{L^{p_{3}}(B_{1/2})}$, $ \|\nabla u_{12}\|_{L^{p_{3}}(B_{1/2})}$ and $ \|\nabla u_{2}\|_{L^{p_{3}}(B_{1/2})}$ as that of (\ref{eq: nonlinear apriori esti-1}) and (\ref{eq: nonlinear apriori esti-2}). Hence
\begin{eqnarray*}
&&\|du\|_{L^{p_{3}}(B_{1/2})}\\
&&\le c\left(\ep+\|f\|_{L^{p}(B_{1})}\right)^{a}\left(\|\na u\|_{M_{\ast}^{4,n-4(1-\al)}(B_{1})}+\|\na^{2}u\|_{M_{\ast}^{2,n-2(2-\al)}(B_{1})}+\|f\|_{L^{p}(B_{1})}\right)
\end{eqnarray*}
for some $a>0$. Similarly, we obtain
\begin{eqnarray*}
\|\na^{2}u\|_{L^{p_{2}}(B_{1/2})}&\leq& c\left(\ep+\|f\|_{L^{p}(B_{1})}\right)^{a}\Big(\|\na u\|_{M_{\ast}^{4,n-4(1-\al)}(B_{1})}\\
&&\qquad\qquad+\|\na^{2}u\|_{M_{\ast}^{2,n-2(2-\al)}(B_{1})}+\|f\|_{L^{p}(B_{1})}\Big)
\end{eqnarray*}
for some $a>0$. Here $c$ and $a$ are two constants that depend on $n, m, p$.


Now we derive the a priori estimate for $\na^{3}u$.
Applying the elliptic regularity theory to the equation (\ref{eq: original equation rewritten}),
we obtain
\[
\|\na^{3}u\|_{L^{p_{1}}(B_{1/2})}\lesssim\|I\|_{L^{p_{1}}(B_{3/4})}
+\|II\|_{L^{p}(B_{3/4})}+\|\nabla u\|_{L^{4}(B_{3/4})}.
\]
By the growth property \eqref{eq: growth of I}, we have
\begin{eqnarray}
\|I\|_{L^{p_1}(B_{3/4})}&\lesssim&\left(\|\na^{2}u\|_{L^{\frac{n}{2}}(B_{3/4})}
+\|\na u\|_{L^n(B_{3/4})}^{2}\right)
\|\na u\|_{L^{p_{3}}(B_{3/4})}\nonumber\\
&\lesssim&(\ep+\|f\|_{L^p(B_1)})^{a}\|\na u\|_{L^{p_{3}}(B_{3/4})}.\nonumber
\end{eqnarray}
Thus we obtain
\[
\|I\|_{L^{p_{1}}(B_{1/2})}\lesssim(\ep+\|f\|_{L^p(B_1)})^{a}\left(\|\na u\|_{M^{4,n-4}(B_{1})}+\|\na^{2}u\|_{M^{2,n-4}(B_{1})}
+\|f\|_{L^{p}(B_{1})}\right).
\]
Similarly, we obtain
\[
\|II\|_{L^{p}(B_{1/2})}\lesssim(\ep+\|f\|_{L^p(B_1)})^{a}\left(\|\na u\|_{M^{4,n-4}(B_1)}+\|\na^{2}u\|_{M^{2,n-4}(B_1)}
+\|f\|_{L^{p}(B_1)}\right).
\]
In conclusion, we deduce
\begin{eqnarray*}
\|\na^{3}u\|_{L^{p_{1}}(B_{1/2})}&\lesssim&
(\ep+\|f\|_{L^p(B_1)})^{a}\Big(\|\na u\|_{M^{4,n-4}(B_{1})}+\|\na^{2}u\|_{M^{2,n-4}(B_{1})}
+\|f\|_{L^{p}(B_{1})}\Big)\\
&&\qquad+\|\nabla u\|_{L^{4}(B_{1})}.
\end{eqnarray*}
This finishes the proof of Theorem  \ref{thm: main result} in
the case $n/4<p<n/3$.

\subsection{Case 2: $n/3\le p<\wq$}

In the remaining case $n/3\le p<\wq$, the result follows by an induction
argument and a trivial iteration.

Since $n/3\le p$, it follows that $f\in L^{q}(B_{1})$ for any $q<n/3$.
Choose $\ep=\ep_{q}$ sufficiently small such that we can apply the
previous result to obtain $u\in W_{\loc}^{3,q}$, which then implies
that
\begin{eqnarray*}
\na u\in L_{\loc}^{s},\na^{2}u\in L_{\loc}^{n-\de} &  & \forall\;s<\wq,0<\de\ll 1.
\end{eqnarray*}
Write the equation as
\[
\De^{2}u={\rm div}\left(I\right)+II+f.
\]
As a result, $I\in L_{\loc}^{n-\de}$  and
$II\in L_{\loc}^{n/2-\de}$ for any $\de>0$ small. Let
\[
\De^{2}u_{1}={\rm div}I,\quad\De^{2}u_{2}=II,\quad\De^{2}u_{3}=f.
\]
We find that $u_{1}\in W_{\loc}^{3,n-\de}$, $u_{2}\in W_{\loc}^{4,\frac{n}{2}-\de}$
and $u_{3}\in W_{\loc}^{4,p}$.

Case 2.1.  If $n/3\le p<n/2$, then for $\de>0$ sufficiently small $W^{3,n-\de}\subset  W^{3,\frac{np}{n-p}}$. Hence in this case $u_{1}+u_{2}+u_{3}\in W_{\loc}^{3,\frac{np}{n-p}}$.
So
\[
n/3\le p<n/2\Rightarrow u\in W_{\loc}^{3,\frac{np}{n-p}}.
\]

Case 2.2. If $p\ge n/2$, then $f\in L^{q}(B_{1})$ for any $q<n/2$. Apply
the above result gives $u\in W_{\loc}^{3,n-\de},$which in turn implies that
$\na u\in L^{\wq}$ and $\na^{2}u\in L_{\loc}^{s}$ for any $s<\wq$. Hence
$I,II\in L_{\loc}^{s}$ for any $s<\wq$. This then gives $u_{1}\in W_{\loc}^{3,s}$,
$u_{2}\in W_{\loc}^{4,s}$ for any $s<\wq$. However, recall that $u_{3}\in W_{\loc}^{4,p}$.
So we can conclude that
\[
\begin{cases}
u\in W_{\loc}^{3,\frac{np}{n-p}} & \text{if }n/2\le p<n,\\
u\in W_{\loc}^{3,s} & \text{for any }s<\wq\text{ if }p\ge n.
\end{cases}
\]

The a priori estimates in this case can be derived similarly and thus omitted. The proof of Theorem  \ref{thm: main result} is complete.

\subsection{Case 3: $1<p\le n/4$}

\begin{proof}[Proof of Theorem \ref{thm:p small}]
The proof of this theorem is almost the same as that of Theorem \ref{thm: main result},
only with minor modifications in the arguments. First note that our
Morrey estimate holds as well. So we can iterate. By the assumption
of $f$, we have
\[
I_{2}(f)\in L^{\eta q}\cap M_{\ast}^{\eta,n-\eta(2-\al)}(\R^{n}),
\]
Remark that $2<\eta<\chi$. This term determines how much regularity
we can gain in the end.

If $\eta q\le\chi$, the iteration stops at the first step, and gives
\[
\na^{2}u\in L^{q\eta}\cap M_{\ast}^{\eta,n-\eta(2-\al)}(B_{1/2}).
\]
In case $\eta q>\chi$, using the same iteration method with slightly
modification, we can obtain the same result. As a result, it follows
from the potential theory that
\[
\na u\in L^{q\eta\chi}\cap M_{\ast}^{\eta\chi,n-\eta\chi(1-\al)}(B_{1/2}).
\]
We leave the details to interested
readers.
\end{proof}

\subsection{Proofs of other results}

\begin{proof}[Proof of Corollary \ref{thm: main result 2}]
	 The proof is standard and omitted here;
	 see for instance \cite[Proposition 6.2]{Guo-Xiang-Zheng-2021-CV}.
\end{proof}

\begin{proof}[Proof of Corollary \ref{thm: Weak compactness}] It follows easily from a contradiction argument; see for instance \cite[Proof of Corollary 1.5]{Guo-Xiang-Zheng-2021-CV}.	
\end{proof}
	
\appendix

\section{Some apriori estimates concerning gauge transform}

Lemma \ref{lem: Gauge transform-general-version} can be proved following
the strategy of Rivi\`ere \cite{Riviere-2007} and Rivi\`ere-Struwe \cite{Riviere-Struve-2008}. We sketch the
proof for readers' convenience. Also, for future applications, we
will prove a  slightly more general  result than that of Lemma \ref{lem: Gauge transform-general-version}.

Let $D\subset\R^{n}$ be a bounded Lipschitz domain, $1<p<\wq$, $1\le q\le\wq$
and $0\le s\le n$. We slightly extend the notion of Morrey spaces.
Say that a function $f$ belongs to the Lorentz-Morrey space $LM^{p,q,s}(D)$,
if $f$ belongs to the Lorentz space $L^{p,q}(D)$, and if
\[
\|f\|_{LM^{p,q,s}(D)}\equiv\sup_{x\in D,0<r<d_{D}}\left(r^{-s/p}\|f\|_{L^{p,q}(B_{r}(x)\cap D)}\right)<\wq,
\]
where $d_{D}$ is the diameter of $D$. Note that
\begin{eqnarray*}
LM^{p,p,s}(D)=M^{p,s}(D) & \text{and} & LM^{p,\wq,s}(D)=M_{\ast}^{p,s}(D).
\end{eqnarray*}
When $s=0$, we get the usual Lorentz space, i.e., $LM^{p,q,0}(D)=L^{p,q}(D)$.
When $0<s\le n$ and $D$ is a bounded domain, we have the continuous
embedding $LM^{p,q,s}(D)\subset L^{p,q}(D)$. Moreover,
\begin{equation}
\|f\|_{L^{p,q}(D)}\le d_{D}^{s/p}\|f\|_{LM^{p,q,s}(D)}\label{eq: a simple embedding}
\end{equation}

\begin{lemma} \label{lem: Apriori Lorentz-Morrey estimate}Let $D\subset\R^{n}$
be a bounded Lipschitz domain, $1<p<\wq$, $1\le q\le\wq$ and $0\le s<n$.
Then, there exists a constant $C>0$ depending only on $D$, $p,q,s$
such that whenever $u\in W_{0}^{1,p}(D)$ is the solution of the equation
\begin{eqnarray*}
-\De u={\rm div}f &  & \text{in }D,
\end{eqnarray*}
for some $f\in LM^{p,q,s}(D,\R^{n})$, then $\na u\in LM^{p,q,s}(D)$.
Moreover,
\[
\|\na u\|_{LM^{p,q,s}(D)}\le C\|f\|_{LM^{p,q,s}(D)}.
\]
\end{lemma}
\begin{proof}
Wen $s=0$ and $q=p$, the result is well known. So the result follows
from a standard interpolation arguments in the case $s=0$ and $1\le q\le\wq$.
In the below we suppose $0<s<n$.

Let $x_{0}\in D$ be an arbitrary point in $D$ and $r>0$. Denote
$D_{r}=D\cap B_{r}(x_{0})$. Let $v$ be the harmonic function in
$D_{r}$ with Dirichlet boundary value $u$. Then, the function $w=u-v$
solves
\[
\begin{cases}
-\De w={\rm div}f & \text{in }D_{r},\\
w=0 & \text{on }\pa D_{r}.
\end{cases}
\]
So apply the result for $s=0$, we obtain
\[
\|\na w\|_{L^{p,q}(D_{r})}\le C\|f\|_{L^{p,q}(D_{r})}.
\]
By the assumption, we find
\[
\|\na w\|_{L^{p,q}(D_{r})}\le C\|f\|_{LM^{p,q,s}(D)}r^{s/p}.
\]
On the other hand, for any $0<\rho<r$,
\[
\|\na v\|_{L^{p,q}(D_{\rho})}^{p}\le C\left(\frac{\rho}{r}\right)^{n}\|\na v\|_{L^{p,q}(D_{r})}^{p}.
\]

Thus, for any $0<\rho<r$, using a simple triangle inequality gives
\[
\|\na u\|_{L^{p,q}(D_{\rho})}^{p}\le C\left(\frac{\rho}{r}\right)^{n}\|\na u\|_{L^{p,q}(D_{r})}^{p}+C\|\na w\|_{L^{p,q}(D_{r})},
\]
from which we deduce that
\[
\|\na u\|_{L^{p,q}(D_{\rho})}^{p}\le C\left(\frac{\rho}{r}\right)^{n}\|\na u\|_{L^{p,q}(D_{r})}^{p}+C\|f\|_{LM^{p,q,s}(D)}^{p}r^{s}.
\]
Therefore, using an elementary lemma, we derive, for any $0<\rho<d_{D}$,
\[
\|\na u\|_{L^{p,q}(D_{\rho})}^{p}\le C\rho^{s}\left(\frac{1}{d_{D}^{s}}\|\na u\|_{L^{p,q}(D)}^{p}+\|f\|_{LM^{p,q,s}(D)}^{p}\right).
\]
Since $x_{0}$ is arbitrary, this is equivalent to
\[
\|\na u\|_{LM^{p,q,s}(D)}\le C\left(\|\na u\|_{L^{p,q}(D)}+\|f\|_{LM^{p,q,s}(D)}\right).
\]
Finally, note that by the result for $s=0$, we have
\[
\|\na u\|_{L^{p,q}(D)}\le C\|f\|_{L^{p,q}(D)}\le C\|f\|_{LM^{p,q,s}(D)}.
\]
The second inequality follows from (\ref{eq: a simple embedding}).
Hence, we conclude from the above two estimates that the desired estimate
holds. The proof is finished.
\end{proof}
Next we consider the following special Poisson equation.

\begin{lemma} \label{lem: apriori BMO-Morrey estimate}Let $D\subset\R^{n}$
be a bounded Lipschitz domain, $1<p<\wq$, $1\le q\le\wq$ and $0\le s<n$.
Then, there exists a constant $C>0$ depending only on $D$, $p,q,s$
such that whenever $u\in W_{0}^{1,p}(D,\wedge^{n-2}\R^{n})$ is the
solution of the equation
\begin{eqnarray*}
-\De u=\ast(dP^{-1}\wedge dP) &  & \text{in }D,
\end{eqnarray*}
for some function $P\in BMO(D)$ with $dP\in LM^{p,q,s}(D)$, then
$du\in LM^{p,q,s}(D)$. Moreover,
\[
\|du\|_{LM^{p,q,s}(D)}\le C\|P\|_{BMO(D)}\|dP^{-1}\|_{LM^{p,q,s}(D)}.
\]
\end{lemma}
\begin{proof}
(1) Suppose $q=p$ and $s=0$, i.e., $P\in BMO(D)$ and $dP\in L^{p}(D)$.
Let $F=|du|^{p-2}du\in L^{p^{\prime}}(D,\wedge^{n-1}\R^{n})$. Hodge
decomposition gives $\psi\in W_{T}^{1,p^{\prime}}(D,\wedge^{n-2}\R^{n})$,
$\be\in W_{N}^{1,p^{\prime}}(D,\wedge^{n-2}\R^{n})$ and an $n-2$
harmonic form $h\in{\cal H}^{n-2}(D,\R^{n})$ such that
$F=d\psi+d^{\ast}\be+h$ and
\[
\|d\psi\|_{p^{\prime}}+\|h\|_{p^{\prime}}\le C\|F\|_{p^{\prime}}=C\|du\|_{p}^{p-1}.
\]
Then
\[
\int_{D}|du|^{p}=\int_{D}du\cdot(d\psi+d^{\ast}\be+h)=\int_{D}du\cdot d\psi.
\]
Here in last equality we used the boundary condition $u=0$ on $\pa D$.
Therefore, we obtain
\[
\int_{D}|du|^{p}=\int_{D}dP^{-1}\wedge dP\wedge\psi=\int_{D}dP^{-1}(P-P_{D})\wedge d\psi,
\]
where $P_{D}=\fint_{D}P$. Since $dP^{-1}\wedge d\psi$ belongs to
Hardy space, we obtain
\[
\int_{D}|du|^{p}\le C\|P\|_{BMO(D)}\|dP^{-1}\|_{L^{p}(D)}\|d\psi\|_{L^{p^{\prime}}(D)}.
\]
This gives
\[
\|du\|_{p}\le C\|P\|_{BMO(D)}\|dP^{-1}\|_{L^{p}(D)}.
\]

(2) In the case $s=0$ and $1\le q\le\wq$, we use the usual interpolation
argument to obtain
\[
\|du\|_{L^{p,q}(D)}\le C\|P\|_{BMO(D)}\|dP^{-1}\|_{L^{p,q}(D)}.
\]

(3) Now suppose $0<s<n$. Use the same arguments as in the Lemma \ref{lem: Apriori Lorentz-Morrey estimate}.
For any $x_{0}\in D$ and $r>0$, denote $D_{r}=D\cap B_{r}(x_{0})$.
Let $v$ be the harmonic function in $D_{r}$ with Dirichlet boundary
value $u$. Then, the function $w=u-v$ solves
\[
\begin{cases}
-\De w=\ast(dP^{-1}\wedge dP) & \text{in }D_{r},\\
w=0 & \text{on }\pa D_{r}.
\end{cases}
\]
Thus using the result in the second step yields
\[
\|dw\|_{L^{p,q}(D_{r})}\le C\|P\|_{BMO(D)}\|dP^{-1}\|_{L^{p,q}(D_{r})}.
\]
It follows
\[
r^{-s/p}\|dw\|_{L^{p,q}(D_{r})}\le C\|P\|_{BMO(D)}\|dP^{-1}\|_{LM^{p,q,s}(D)}.
\]
On the other hand, for any $0<\rho<r$,
\[
\|dv\|_{L^{p,q}(D_{\rho})}^{p}\le C\left(\frac{\rho}{r}\right)^{n}\|dv\|_{L^{p,q}(D_{r})}^{p}.
\]

Therefore, a similar argument as in the previous lemma gives, for
any $x_{0}\in D$ and $0<\rho<d_{D}$,
\[
\|du\|_{L^{p,q}(D_{\rho})}\le C\rho^{s/p}\left(\|du\|_{L^{p,q}(D)}+\|P\|_{BMO(D)}\|dP^{-1}\|_{LM^{p,q,s}(D)}\right).
\]
Since $\|dP^{-1}\|_{L^{p,q}(D)}\le C\|dP^{-1}\|_{LM^{p,q,s}(D)}$,
using the result in the second step together with the above estimate,
we deduce the desired estimate. The proof is complete.
\end{proof}
Based on the above two lemmata, we can prove Lemma \ref{lem: Gauge transform-general-version}.
We prove a slightly more general result here.

\begin{lemma}\label{lem: higher order apriori estimate for gauge}
There exist $\de>0$ and $C>0$ with the following property: Suppose
that $\Om\in LM^{p,q,s}(B_{1/2})$ for some $1<p<\wq$, $1\le q\le\wq$
and $0<s<n$ such that there exist $P,\xi\in LM^{p,q,s}(B_{1/2})$
satisfying the equation \eqref{eq: gauge transform equation} of Lemma
\ref{lem: Gauge transform} on $B_{1/2}$, and
\[
\|dP\|_{M^{4,n-4}(B_{1/2})}+\|d\xi\|_{M^{4,n-4}(B_{1/2})}\le\de,
\]
then there hold
\[
\|dP\|_{LM^{p,q,s}(B_{1/2})}+\|d\xi\|_{LM^{p,q,s}(B_{1/2})}\le C\|\Om\|_{LM^{p,q,s}(B_{1/2})}.
\]

If, in addition, $\na\Om\in M^{\frac{p}{2},\frac{n-p+s}{2}}(B_{1/2})$,
then $\na^{2}P,\na^{2}\xi\in M^{\frac{p}{2},\frac{n-p+s}{2}}(B_{1/2})$,
and
\[
\left\|\na^{2}P\right\|_{M^{\frac{p}{2},\frac{n-p+s}{2}}(B_{\frac{1}{2}})}+\left\|\na^{2}\xi\right\|_{M^{\frac{p}{2},\frac{n-p+s}{2}}(B_{\frac{1}{2}})}\le C\left(\left\|\na\Om\right\|_{M^{\frac{p}{2},\frac{n-p+s}{2}}(B_{\frac{1}{2}})}+\left\|\Om\right\|_{LM^{p,q,s}(B_{\frac{1}{2}})}\right).
\]
In particular, \eqref{eq: improved gauge regularity 1} and \eqref{eq: improved gauge regularity 2}
holds under the assumption $\Om\in M_{\ast}^{4,n-4+4\al}(B_{1/2})$
and $\na\Om\in M_{\ast}^{2,n-4+2\al}(B_{1/2})$.\end{lemma}
\begin{proof}
By equation (\ref{eq: gauge transform equation}),
\[
\begin{cases}
\De\xi=\ast dP^{-1}\wedge dP+\ast d(P^{-1}\Om P) & \text{in }B_{1/2},\\
\xi=0 & \text{on }B_{1/2}.
\end{cases}
\]
Let $\xi_{1}$ be the solution of
\begin{eqnarray*}
\begin{cases}
\De\xi_{1}=\ast dP^{-1}\wedge dP & \text{in }B_{1/2},\\
\xi=0 & \text{on }B_{1/2},
\end{cases} & \text{and} & \begin{cases}
\De\xi_{2}=\ast d(P^{-1}\Om P) & \text{in }B_{1/2},\\
\xi=0 & \text{on }B_{1/2}.
\end{cases}
\end{eqnarray*}
Applying Lemma \ref{lem: Apriori Lorentz-Morrey estimate} to $\xi_{2}$
and Lemma \ref{lem: apriori BMO-Morrey estimate} to $\xi_{1}$, we
deduce
\[
\|d\xi_{1}\|_{LM^{p,q,s}(B_{1/2})}\le C\de\|dP\|_{LM^{p,q,s}(B_{1/2})}
\]
and
\[
\|d\xi_{2}\|_{LM^{p,q,s}(B_{1/2})}\le C\|\Om\|_{LM^{p,q,s}(B_{1/2})}.
\]
Thus
\[
\|d\xi\|_{LM^{p,q,s}(B_{1/2})}\le C\de\|dP\|_{LM^{p,q,s}(B_{1/2})}+\|\Om\|_{LM^{p,q,s}(B_{1/2})}.
\]
Directly from equation (\ref{eq: gauge transform equation}), we have
\[
\|dP\|_{LM^{p,q,s}(B_{1/2})}\le C\|d\xi\|_{LM^{p,q,s}(B_{1/2})}+\|\Om\|_{LM^{p,q,s}(B_{1/2})}.
\]
Combining the above two estimate together with a suitably chosen $\de<1$
small enough, we obtain the first estimate.

The second estimate can be proved by the same method. We omit the
details. The proof is complete.
\end{proof}

\end{document}